\documentclass{article}
\usepackage{latexsym,a4}

\title{ \bf A traditional dealing with a semi-classical limit
and Hopf theorem}
\author{Yanlin Yu}\date{}
\begin{document}
\maketitle

\begin{center}
\parbox{12cm}
{\quad Abstract. \ This paper deals with a semi-classical limit(Theorem 1)
   by using traditional mathematical methods, and shows a Hopf theorem as 
   a corollary. A formal discussing of it may be found in [7]      }
\end{center}

\section{A semi-classical limit theorem}

Let $M$ be a compact, closed Riemannian manifold of dim $n$, 
and $V$ a vector field without  degenerate zeros on $M$.
Let $\Lambda^\ast (M)$ be the space of differential forms on $M$, and
$$D = d + \delta : \Lambda^{\ast}(M) \rightarrow 
\Lambda^{\ast}(M)$$
be the de Rham-Hodge operator, which is an elliptic operator.
Let us consider a Witten's deformation of 
$d + \delta$.
$$D_t = (d+\delta) + t[V^\ast \wedge + i(V)]: \Lambda^{\ast} (M) 
\rightarrow \Lambda^{\ast}(M),$$
where $V^\ast$ is a 1-form dual to the vector field $V$, and $V^\ast \wedge$
 means the exterior product by $V^\ast$, while $i(V)$ the interior 
product by $V$. Let
$$\Box_t = D^2_t:\Lambda^{\ast}(M) \rightarrow \Lambda^{\ast}(M),$$
and $e^{-\tau \Box_t}$ 
be the solution operator of the heat operator
$\frac {\partial} {\partial \tau} + \Box_{t}.$
 It is well known that $e^{-\tau \Box_t}$ is an integral operator, i.e.
 there exists a unique family of linear maps
$$G(\tau, q, p, t): \Lambda^\ast_{p}(M) \rightarrow \Lambda^\ast_{q}(M)$$
such that 
$$(e^{-\tau \Box_t} \phi)(q) =
\displaystyle\int_M G(\tau, q, p, t) \phi (p) dp,
\quad\quad \forall \phi.$$
Such a family of $G(\tau, q, p, t)$ 
 is  called a fundamental solution of the heat operator 
$\frac{\partial}{\partial\tau} + \Box_t$. The fundamental solution can be 
determined by the following equations 
$$\left\{ \begin{array}{l}
 \left(\frac{\partial}{ \partial  \tau } +\Box _t \right) G(\tau, 
q, p, t) = 0 \\[3mm] 
\displaystyle\lim _{\tau \to 0 } \displaystyle\int _M G(\tau,q,p,t) 
\phi (p) d p=\phi (q), \quad \forall \phi,
\end{array}  \right.    $$
where $\Box_t$ acts on the indeterminate $q$.
 If $\Box_{t}$ is thought as a 
deformation of a physical system, Witten ([3]) had considered a limit 
situation of $\Box_{t}$ as $t\rightarrow \infty$, and called it  a 
``semi-classical limit''. By using this consideration he and [2] gave a proof 
of Morse inequalities. Afterwards Bismut([1]) considered a double limit 
$$\displaystyle\lim_{s \rightarrow \infty}
\displaystyle\lim_{\tau \rightarrow 0}\hbox{str}e^{-\tau \Box_{\frac{s}{\tau}}}
$$
to give another proof of Morse inequalities, where str means a super
trace we will explain later. The first limit 
$\displaystyle\lim_{\tau \rightarrow 0}$  actually means that 
$$''\tau \rightarrow 0  \quad \hbox{and}\quad t \rightarrow \infty \quad 
\hbox{and}
\quad  \tau t \quad \hbox{is kept as a constant} .''$$
We call this limit a semi-classical limit too, and denote it by $s-\lim$.  
As the semi-classical limit
$$(s-lim) \hbox{str}e^{-\tau \Box_{\frac{s}{\tau}}}=
(s-lim)\displaystyle \int_{M} \hbox{str} \,G(\tau,p,p,t)dp $$
is concerned, how to understand the fundamental solution $G(\tau,q,p,t)$
is a very serious thing. We introduce $\Phi_{0}(\tau,t,p)$ as follows
in order to replace $G(\tau,p,p,t)$ when the semi-classical limit is
evaluated.

For any $p\in M$,choosing an orthonormal frame 
$\{E_1(p), \cdots , E_n(p) \}$ at $p$,
thus the vector $V(p)$ can be expressed as 
$$V(p) = \displaystyle\sum_{i} v_i(p) E_i(p).$$ 
Define $v_{ij}(p)$ by 
$$\bigtriangledown_{E_j(p)} V = \displaystyle\sum_{i} v_{ij}(p) E_{i}(p),$$
where $\bigtriangledown$ is the Levi-Civita connection.
In general the matrix $(v_{ij}(P))$ is not symmetric, we denote it by $A(p)$.
Let $A(p)^{*}$ be the transpose of $A(p)$,
$$\theta = \theta(\tau,t,p)=2\tau t \sqrt{A(p)A(p)^{*}}.$$
Define a linear map
$$\phi_{0} (\tau,t,p):\Lambda^{\ast}_{p}(M) \rightarrow\Lambda^{\ast}_{p}(M)$$
by
$$\begin{array}{rl}
\phi_{0}(\tau,t,p)&=\frac{1}{\sqrt{4 \pi \tau}^n}\sqrt 
{\det\left(\frac{\theta}{\sinh\theta}\right) }\cr
&\cdot \exp \left\{-2\tau t^{2}(v_{1}(p),\cdots,v_{n}(p))
 \frac{\cosh\theta -1}{\theta \sinh\theta}
\left(\begin{array}{c}
v_{1}(p)\\ \vdots\\ v_{n}(p)
\end{array}\right)
+\tau t\sum v_{ij}(p)E^{+}_{i}E^{-}_{j}\right\},
\end{array}$$
 where
 $$E^{\pm}_{j}= \omega_{j}(p) \pm i(E_{j}(p)) :
  \Lambda^{*}_{p}(M) \rightarrow \Lambda^{*}_{p}(M),$$
  $\{\omega_{1}(p),\cdot, \omega_{n}(p)\}$ 
is the coframe dual to
    $\{E_{1}(p),\cdot, E_{n}(p)\}.$ It is easy to see that 
$\phi_{0}$ does not depend on the choice of 
$\{ E_{1}(p),\cdots,E_{n}(p)\}.$
In this paper we will prove the following theorem

{\bf Theorem 1} \ Let $M$ be a compact closed Riemannian manifold, $V$ a 
vector field without degenerated zeros. Then   
 $$ (s-\lim) \displaystyle \int_{M}|G(\tau,p,p,t)- \phi_{0}(\tau,t,p)|dp =0,$$
where we used the norm of a linear map,which is defined as usually, i.e.
$$|\psi|=\sqrt{\hbox{tr}(\psi \psi^{\ast})}.$$
In order to prove theorem 1, we need to introduce a parametrix out of 
considerations of harmonic oscillators in {\S 2}, and by using Lemma A and 
Lemma B in {\S 4} we can compare the parametrix with $G(\tau,q,p,t)$. 
The proof of theorem 1 will be finished in {\S 6}. In {\S 7} we prove the
Hopf theorem. The appendix, which shows an independent
interest,  is needed when we prove lemma A.

\section{ Harmonic oscillators}
Let us recall Mehler formula first. It is
$$\begin{array}{rl}
{\cal M}(\tau , y, x, b)=&\frac{1}{\sqrt{4 \pi \tau}}
\sqrt{\left(\frac{\theta}{\sinh\theta}\right)}\\[3mm]
&\cdot \exp \left\{- \frac{1}{4 \tau} \frac{\theta}
{\sinh\theta} [\cosh\theta\cdot (x^{2} + y^{2})-2xy] \right\},
\end{array}$$
where $\theta=2b\tau.$ The Mehler formula satisfies

$$\left\{\begin{array}{l} 
\left[\frac{\partial}{\partial \tau}- \frac{\partial^2}{\partial y^2}
   + b^2 y^2\right]{\cal M}(\tau,y,x,b)=0 \\[3mm]
\displaystyle\lim_{\tau \rightarrow 0} \displaystyle\int_R
{\cal M} (\tau,y,x,b) f(y) dz = f(x). 
\end{array}\right. $$
From the above equations, it is easy to see that if
$${\cal M}_0 (\tau,y,x,b) = {\cal M} (\tau,y+x,x,b),$$
then
$$\left\{\begin{array}{l}
\left[\frac{\partial}{\partial \tau} -\frac{\partial^2}{\partial y^2}
   + b^2 (y+x)^2\right]{\cal M}_0(\tau,y,x,b)=0 \\[3mm]
\displaystyle\lim_{\tau \rightarrow 0} \displaystyle\int_R
{\cal M}_0 (\tau,y,x,b) f(y) dy = f(0). 
\end{array}\right.$$
Again, if
$$\Phi(\tau,y,a,b)= {\cal M}_0 (\tau,y,\frac{a}{b},b),$$ then holds
$$\left\{\begin{array}{l}
\left[\frac{\partial}{\partial \tau} -\frac{\partial^2}{\partial y^2}
   + (a+by)^2)\right]\Phi (\tau,y,a,b)=0 \\[3mm]
\displaystyle\lim_{\tau \rightarrow 0} \displaystyle\int_R
\Phi (\tau,y,a,b) f(y) dy = f(0). \end{array}\right.$$
Let us check wheather $\Phi$ has a singularity at $b=0$.Due to
$$\begin{array}{rl}
\Phi(\tau,y,&a,b)={\cal M}_0 (\tau,y,\frac{a}{b},b)=
    {\cal M}_0 (\tau,y+\frac{a}{b},\frac{a}{b},b)     \\[3mm]
=&\frac{1}{\sqrt{4 \pi \tau}}
  \sqrt{\frac{\theta}{\sinh\theta}}
  \hbox{exp}
   \{-\frac{1}{4\tau} 
      \left[\frac{\theta\cosh\theta}{\sinh\theta}
  (\frac{a^2}{b^2}+(y+\frac{a}{b})^2) -2
   \frac{\theta }{\sinh\theta}\frac{a}{b}
   (y+\frac{a}{b}) 
   \right]\}\\
=&\frac{1}{\sqrt{4 \pi \tau}}
  \sqrt{\frac{\theta}{\sinh\theta}}
  \hbox{exp}
   \{-\frac{1}{4\tau} 
      \left[\frac{\theta\cosh\theta}{\sinh\theta}
  y^2
  +4a\tau \frac{\cosh\theta -1}{\sinh\theta} y
  +8a^2 \tau^2 \frac{\cosh\theta -1}{\theta\sinh\theta}
 \right]\},     
\end{array}$$
there is no singularity at $b=0$!
From the above formula of $\Phi(\tau,y,a,b)$ it follows that
$$\begin{array}{rl}
\Phi(\tau,y,a,b)=
&\frac{1}{\sqrt{4 \pi \tau}}
  \sqrt{\frac{\theta}
             {\sinh\theta}
       }
  \hbox{exp}\{-\frac{1}{8\tau} 
               \frac{\theta(\cosh\theta+1)}
                    {\sinh\theta}y^2 
\\&\,\,\,\,\,\,\,\,
             -\frac{1}{4\tau}
               \left[a\tau 
                   \sqrt{8\frac{\cosh\theta-1}
                               {\theta \sinh\theta}
                        }
          +y\sqrt{\frac{\theta(\cosh\theta-1)}
                       {2\sinh\theta}
                 }
          \right]^2\} \\
\leq&\frac{1}{\sqrt{4 \pi \tau}}
  \sqrt{\frac{\theta}
             {\sinh\theta}
       }
  \hbox{exp}\{-\frac{1}{8\tau} 
               \frac{\theta(\cosh\theta+1)}
                    {\sinh\theta}y^2 \} \\
\Phi(\tau,y,a,b)=            
&\frac{1}{\sqrt{4 \pi \tau}}
  \sqrt{\frac{\theta}{\sinh\theta}
       }
  \hbox{exp}\{-\tau a^2 \frac{\sinh\theta}
                             {\theta\cosh\theta}
  -\frac{1}{4\tau}\left[\sqrt{\frac{\theta\cosh\theta} 
                              {\sinh\theta}} y       
      +2a\tau\frac{\cosh\theta-1}
                  {\sqrt{\theta\cosh\theta \cdot\sinh\theta}
                  }\right]^{2}
           \}\\
\leq&\frac{1}{\sqrt{4 \pi \tau}}
  \sqrt{\frac{\theta}{\sinh\theta}
       }
  \hbox{exp}\{-\tau a^2 \frac{\sinh\theta}
                             {\theta\cosh\theta}\}.   
\end{array}$$

{\bf Definition 2} For a vertor $a=(a_1,\cdots,a_n) \in R^n$, and an 
$n\times n$
matrix $B$, let $\Theta$ and $\Theta^\#$ be arithmetic square roots of
$4\tau^2 BB^{\ast}$ and $4\tau^2 B^{\ast}B$ respectively, where $B^{\ast}$ is 
the transpose of $B$.
 We define
$$\begin{array}{rl}
\Phi(\tau,Y,&a,B)=\frac{1}{\sqrt{4 \pi \tau}^{n}}         
  \sqrt{\det\left(\frac{\Theta}{\sinh\Theta}\right)} \\[3mm]
 & \hbox{exp}\{-\frac{1}{4\tau} Y \frac{\Theta\cosh\Theta}
{\sinh\Theta}Y^{\ast}
  -2\tau Y \frac{\cosh\Theta -1}{\Theta\sinh\Theta} Ba^{\ast}
  -2 \tau a \frac{\cosh\Theta^\# -1}{\Theta^\#\sinh\Theta^\#} a^{\ast}
 \},     
\end{array}$$
where $Y=(y_1,\cdots,y_n), Y^{\ast}$ and $a^{\ast}$ are the transposes of 
$Y$ and $a$ respectively.

{\bf Proposition 3}\quad  $\Phi (\tau,Y,a,B)$ satisfies
$$\left\{\begin{array}{l}\left[\frac{\partial}{\partial \tau} -
                \sum^{n}_{i=1}\frac{\partial^2}{\partial y^2}
               + (a+YB)(a+YB)^{\ast}\right]\Phi (\tau,Y,a,B)=0 \\[3mm]
         \displaystyle\lim_{\tau \rightarrow 0} \displaystyle\int_R
              \Phi (\tau,Y,a,B) f(Y) dY = f(0). \end{array}\right.$$

{\bf Proof}\ Sometime we denote an $(i,j)-$element of a matrix $C$ by $C_{ij}$,
and $\cosh\Theta$ by $\cosh$, and $\cosh\Theta^\#$ by $\cosh^\#$.
 Note that $\frac{\Theta}{\tau}$ does not depend on $\tau$, 
hence $|frac{\partial}{\partial \tau}(\frac{\Theta}{\tau})=0.$ This fact  helps 
the following computations. First
$$\begin{array}{rl}
\frac{\partial}{\partial \tau}\Phi =& \Phi \left\{
           { \sqrt{\det\left( \frac{\Theta}{\tau \sinh}\right) }}
^{-1}
            \frac{\partial}{\partial \tau}
            \sqrt{\det\left(    \frac{\Theta}{\tau \sinh}\right)}
        -\frac{1}{4}Y \frac{\Theta}{\tau} \frac{\partial}{\partial \tau}
        \left(\frac{\cosh} {\sinh}    \right) Y  \right.    \\  
          & \left. -2Y \frac{\tau}{\Theta} \frac{\partial}{\partial \tau}
           \left(   \frac{\cosh-1}{\sinh}   \right)
        Ba^{\ast}    
  -2a \frac{\tau}{\Theta^\#} \frac{\partial}{\partial \tau}
   \left( \frac{\cosh^\#-1}{\sinh^\#} \right)
   a^{\ast}  \right\},
\end{array}$$
then
$$\begin{array}{rl}
 {\sqrt{\det\left(    \frac{\Theta}{\tau \sinh}\right)}}^{-1}
&\frac{\partial}{\partial \tau}
  \sqrt{\det\left(    \frac{\Theta}{\tau \sinh}\right)}=
\frac{1}{2}\left(\det\left(\frac{\Theta}{\tau \sinh}\right)
         \right)^{-1}   \frac{\partial}{\partial \tau}
         \left( \det\left(\frac{\Theta}{\tau \sinh}\right)\right)\\
=&       \frac{1}{2}\hbox{tr}\left(
         \left(\frac{\Theta}{\tau \sinh}\right)^{-1}
         \frac{\partial}{\partial \tau}
         \left(\frac{\Theta}{\tau \sinh\Theta}\right)
         \right)    
=\frac{1}{2}\hbox{tr}\left(\frac{\tau \sinh}{\Theta}
                         \frac{\Theta}{\tau}
                         \left(
                         -\frac{\cosh}{\sinh^{2}}
                         \frac{\Theta}{\tau}
                         \right)
                     \right)
                     \\
              =&-\frac{1}{2}\hbox{tr}\left(
                 \frac{\Theta \cosh}{\tau \sinh}
                 \right).
\end{array}$$
And then by
$$\frac{\Theta}{\tau} \frac{\partial}{\partial \tau}
   \left( \frac{\cosh}{\sinh}\right)
=-\frac{\Theta^2}{\tau^2 \sinh^2},
   \quad \quad \quad 
\frac{\tau}{\Theta} \frac{\partial}{\partial \tau}
   \left( \frac{\cosh-1}{\sinh}\right)
    =\frac{\cosh-1}{\sinh^2},$$
 we get
 $$\frac{\partial}{\partial \tau}
   \Phi=\Phi \left\{
-\frac{1}{2 \tau} \hbox{tr}\left(
   \frac{\Theta \cosh}{\sinh}\right)
   +\frac{1}{4 \tau^2} Y \frac{\Theta^2}{\sinh^2} Y^{\ast}
   -2 Y \frac{\cosh-1}{\sinh^2} Ba^{\ast}
   -2a \frac{\cosh^\#-1}{{\sinh^\#}^2}a^{\ast}
   \right\}. $$
   Again,
$$\begin{array}{rl}
\frac{\partial}{\partial y_{i}}\Phi&=\Phi \left\{
-\frac{1}{2 \tau} \left(
   \frac{\Theta \cosh}{\sinh}\right)_{ij}y_{j}
   -2 \tau \left(\frac{\cosh-1}{\Theta \sinh}\right)_{ij} (Ba^{\ast})
_{j}
                                        \right\},        \\[3mm]
\sum_{i}\frac{\partial^{2}}{\partial y^{2}_{i}}&\Phi=\Phi \left\{
-\frac{1}{2 \tau} \hbox{tr}\left(\frac{\Theta \cosh}{\sinh}\right)
+\frac{1}{4 \tau^2} Y 
\left( \frac{\Theta \cosh}{\sinh}\right)^2 
Y^{\ast} \right.\\[3mm]
&+2 Y \frac{\cosh(\cosh-1)}{\sinh^2} Ba^{\ast}
+4\tau^2 aB^{\ast} 
\left( \frac{\cosh-1}{\Theta \sinh}\right)^2  
Ba^{\ast} \left. \right\}.
\end{array}$$
Therefore
$$\begin{array}{rl}
\Phi^{-1}(\frac{\partial}{\partial \tau}\Phi&
-\sum_{i}\frac{\partial^{2}}{\partial y^{2}_{i}}\Phi) =
-\frac{1}{4\tau^2}Y \Theta^2 Y^{\ast}
-2 YBa^{\ast}
-2 a\frac{\cosh^\# -1}{\sinh^{\#2}}a^{\ast}
\\&\,\,\,\,\,\,\,\,
-4 \tau^2 a B^{\ast} \left(\frac{\cosh-1}{\Theta\sinh}\right)^2 Ba^{\ast}
    \\[3mm]
&=-\frac{1}{4\tau^2} Y \Theta^2 Y^{\ast}  -2YBa^{\ast} -2a
\frac{\cosh^\#-1}{\sinh^{\#^2}} a^{\ast}
-a\left(\frac{\cosh^\# -1}{\Theta^\#\sinh^{\#}}\right)^2
\Theta^{\#2}a^{\ast}
    \\[3mm]
&=-YBB^{\ast} Y^{\ast} -2 Y B a^{\ast}  -aa^{\ast}   \\[3mm]
&=-(YB+a)(B^{\ast} Y^{\ast}+a^{\ast}).
\end{array}$$
The first half of the proposition is proved. Because the second half is easy,
its proof is omitted.

{\bf Proposition 4} \quad There hold
$$\begin{array}{rl}
(i)&\quad \Phi (\tau,Y,a,B) \leq
\frac{1}{{\sqrt{4 \pi \tau}}^n }
\sqrt{\det\left(\frac{\Theta}{\sinh}\right)}
\hbox{exp}\left\{- \frac{1}{8\tau} Y\frac{\Theta (\cosh+1)}{\sinh} Y^{\ast}
\right\}    \\[3mm]
(ii)& \quad \Phi (\tau,Y,a,B) \leq
\frac{1}{{\sqrt{4 \pi \tau}}^n }
\sqrt{\det\left(\frac{\Theta}{\sinh}\right)}
\hbox{exp}\left\{- \tau a \frac{\sinh}{\Theta \cosh} a^{\ast}
\right\}
\end{array}$$

{\bf Proof}\quad The proof is similar to the case of $\Phi (\tau,y, a, b),$
which was mentioned before.

{\bf Proposition 5}\quad  Let
$\Phi_0 (\tau,Y,X,B)=\Phi(\tau,Y-X,XB,B),$
then
$$\begin{array}{rl}
(i) 
\Phi_0 (\tau,Y,a,B)=&\frac{1}{{\sqrt{4 \pi \tau}}^{n}}
  \sqrt {\det\left(\frac{\Theta}{\sinh}
                   \right)
         }
 \hbox{exp}  \left\{  -\frac{1}{4\tau} Y
                 \left( \frac{\Theta\cosh}{\sinh}
                 \right) Y^{\ast}                \right.      \\[3mm]
 &
  -\frac{1}{4 \tau} X  \left( \frac{\Theta\cosh}{\sinh}
                    \right) X^{\ast}
 +\frac{1}{2 \tau} Y  \left( \frac{\Theta}{\sinh}
                   \right) X^{\ast}
             \left.\right\} \\
&           {\bf (Mehler} \quad {\bf formula)},\\                       \\[3mm]
(ii) \Phi_0 (\tau,Y,a,B)=&\frac{1}{\sqrt{4 \pi \tau}^{n}}
  \sqrt {\det\left(\frac{\Theta}{\sinh}
                   \right)
         }
 \hbox{exp}\left\{     -\frac{1}{8\tau} (Y+X)
                      \frac{\Theta(\cosh-1)}{\sinh}
                      (Y+X)^{\ast}\right.             \\[3mm]
 &-\frac{1}{8\tau} (Y-X)
                     \frac{\Theta(\cosh+1)}{\sinh}
                     (Y-X)^{\ast}
         \left. \right\}.
 \end{array}$$

{\bf Proof} We only check (i) as follows
$$\begin{array}{rl}
  \hbox{exp}&\{-\frac{1}{4\tau} (Y-X) \frac{\Theta\cosh}{\sinh}(Y-X)^{\ast}
  -2\tau (Y-X) \frac{\cosh -1}{\Theta\sinh} B(XB)^{\ast}
\\ &\,\,\,\,\,\,\,\,
  -2 \tau (XB) \frac{\cosh^\# -1}{\Theta^\#\sinh^\#} (XB)^{\ast}
 \}      \\[3mm]
 &=\hbox{exp}\{-\frac{1}{4\tau} Y  \frac{\Theta\cosh}{\sinh}Y^{\ast}
     -\frac{1}{4\tau} X  \frac{\Theta\cosh}{\sinh}X^{\ast}
\\[3mm]
    & +\frac{1}{2\tau}Y\left(\frac{\Theta\cosh}{\sinh}
     -4 \tau^{2}\frac{\cosh -1}{\Theta\sinh}BB^{\ast}\right)X^{\ast}
     +\frac{1}{2\tau}X\left(4\tau^2  \frac{\cosh -1}{\Theta\sinh} 
     BB^{\ast}\right. 
 \\ &\,\,\,\,\,\,\,\, \left. 
     -4\tau^2 B \frac{\cosh^\# -1}{\Theta^\#\sinh^\#} B^{\ast} \right)X^{\ast}\}
     \\[3mm]
 &=\hbox{exp}\left\{ -\frac{1}{4\tau} Y \frac{\Theta \cosh}{\sinh} Y^{\ast}
      -\frac{1}{4\tau} X \frac{\Theta \cosh}{\sinh} X^{\ast}
      +\frac{1}{2\tau}Y \frac{\Theta}{\sinh}X^{\ast} \right\}.
\end{array}$$
Therefore (i) is true.

\section {Parametrix}
Choose a local orthonormal frames $\{E_{1},\dots,E_{n}\}$ on $M$, let
$\{\omega_{1},\cdots,\omega_{n}\}$ be the coframes dual to
 $\{E_{1},\dots,$  $E_{n}\}$, then 
$$V=v_{i}E_{i}, \qquad V^*=v_{i}\omega_{i},$$
and
$$V^*\wedge +i(V)=v_{i}E^{+}_{i},$$ so
$$D_{t}=(d+\delta) + tv_{i}E_{i}^{+}=E_{i}^{-}\nabla_{E_{i}}
+tv_{i}E_{i}^{+}.$$
Therefore
$$\begin{array}{rl}
\Box_t =D^{2}_{t}&=(d+\delta)^2+\sum_{j,k}
(E_{k}^{-}\nabla_{E_k}tv_{j}E_{j}^{+}+tv_{j}E_{j}^{+}E_{k}^{-}
\nabla_{E_{k}}) + t^2 \sum_{j,k}v_j v_k E_{j}^{+}E_{k}^{+} \\
&=(d+\delta)^2 +\sum_{j,k}tv_{jk}E_{k}^{-}E_{j}^{+}
  +\sum_{j,k}tv_{j}(E_{k}^{-}E_{j}^{+} +E_{j}^{+}E_{k}^{-})
\nabla_{E_{k}}+t^2\sum_{j}v_{j}^{2}  \\
&=(d+\delta)^2 - \sum_{j,k}t v_{jk}E^{+}_{j}E^{-}_{k}
  +t^{2}\sum_{j}v^{2}_{j}.\end{array}$$
(One may see [yu 6] for formulas for the multiplication table of 
$E_{1}^{+},\cdots,E_{n}^{+},$ $E_{1}^{-},\cdots,E_{n}^{-}.$)

Now for $p\in M$ in a neighbourhood of $p$  we choose a normal coordinate 
system $\{y_{1},\cdots,y_{n}\}$  centering at $p$ and an orthonormal moving 
 frame
$\{E_{1},\cdots,E_{n}\}$, which is parallel along geodesics passing through
$p$ and 
$$E_{i}(p)=\frac{\partial}{\partial y_{i}}|_{p}.$$
Of course the coordinates of $p$ is $(0,\cdots,0).$ Suppose the coordinates of $q$ is $(y_{1},\cdots,y_{n})=Y,$ for $q$ near $p.$ $\Box_{t}$ can be written as
$$\Box_t=-\sum_{i}\frac{\partial^{2}}{\partial y_{i}^{2}} 
+t^{2}\sum_{i}(v_{j}(p) +\sum_{k}v_{jk}(p)y_{k})^{2}
-\sum_{j,k}tv_{jk}(p)E_{j}^{+}E_{k}^{-}+\ldots,$$
which suggests the following definition of $H(\tau,q,p,t)$ after comparing 
$\Box_t$
with the first equation in Proprsition 3. $H(\tau,q,p,t)$ is chosen for a good 
approximation of the fumdamental solution $G(\tau,q,p,t).$

{\bf Definition 6}\ Define
$$H(\tau,q,p,t):\Lambda_{p}^{*}(M)\rightarrow\Lambda_{q}^{*}(M)$$ by
$$H(\tau,q,p,t)=\Phi(\tau,Y,(tv_{1}(p),\cdots,tv_{n}(p)),tB(p))
\hbox{exp}\left\{\sum \tau t v_{jk}(p)E_{j}^{+}E_{k}^{-}\right\} 
\cdot\phi(q,p),$$
where $\left(B(p)\right)_{ij}=v_{ji}(p),\phi (q,p)$ is a $C^{\infty}$ fuction, which
equals $1$ in a small neighbourhood of the diagonal of $M\times M$ and $0$
outside a little larger neighbourhood.

The proposition 13 in {\S 5} means that this $H(\tau,q,p,t)$ is a parametrix
of $\frac{\partial}{\partial \tau}+\Box_{t}.$ Let
$$Q(\tau,q,p)=\frac{1}{\sqrt{4\pi \tau}^{n}}\hbox{exp}
\left\{-\frac{\rho(q,p)^2}{4\tau}\right\}:
\Lambda_{p}^{*}(M)\rightarrow\Lambda_{q}^{*}(M),$$
then the proposition 4 in {\S 2} induces the following lemma immediately.

{\bf Lemma 7}\ For $s_{0}>0$ there exist $c_0, c_1 >0$ such that for 
$\tau t <s_0$ we have
$$|H(\tau,q,p,t)|\leq c_0 Q(c_1 \tau, q,p)\hbox{exp}
\left\{-\frac{\tau t^2}{c_1}v(p)^2 \right\}.$$

{\bf Proof}\ Note that
$$|\hbox{exp}\left\{\tau t v_{jk}(p)E_{j}^{+}E_{k}^{-}\right\}|
\leq \hbox{const.},$$ and consider the square root of the product of 
right-hand sides of inequalities in the proposition 4, we get a proof easily.

{\bf Lemma 8}\ For $\tau_0, s_0 >0$ there exist $c_0,c_1$ such that for any 
$\tau,t$ with $0<\tau \leq \tau_0$ and $0<\tau t\leq s_0$ we have
$$|(\frac{\partial}{\partial \tau}+\Box_t)
H(\tau,q,p,t)|\leq c_0 (\sqrt{\tau}t +1)
Q(c_1 \tau, q,p)\hbox{exp}
\left\{-\frac{\tau t^2}{c_1}v(p)^2 \right\}.$$

{\bf Proof}\ First we recall some notations and facts in [5], let
$H^{ij},\Gamma^{k}_{ij},R_{ijkl}$ be defined by 
$$\begin{array}{l}
E_j=H^{ij}\frac{\partial}{\partial y_i},\\
\nabla_{E_{i}}E_{j}=\Gamma^{k}_{ij}E_{k},\\
R_{ijkl}=-<(\nabla_{E_i}\nabla_{E_j}-\nabla_{E_j}\nabla_{E_i}
-\nabla_{[E_i,E_j]})E_k,E_l>,
\end{array}$$ 
and let $\left(H_{ij}\right)$ be the inverse of  $\left(H^{ij}\right)$.
A corollary 8 in [5] claimed a Taylor's expansion 
$$H_{ij}(y)=\delta_{ij}+\frac{1}{6}\sum_{k,l}R_{ijkl}(p)y_k y_l+\ldots.$$
As usual $\Delta_0$ is defined by
$$\Delta_0 =\sum_i (\nabla_{E_i}\nabla_{E_i}-\nabla_{F_i}),$$
where $F_i =\nabla_{E_i}E_i.$  Weizenbock formula reads
$$(d+\delta)^2 =-\Delta_0 +{\tilde R},$$
where 
$${\tilde R}=\frac{1}{8}\sum R_{ijkl}E_{i}^{-}E_{j}^{-}E_{k}^{+}E_{l}^{+}
   +\frac{1}{4}\sum_{i,j}R_{ijij}.$$
Thus
$$\Box_t =-\Delta_0  - \sum_{jk}t v_{jk}E^{+}_{j}E^{-}_{k}
  +t^{2}\sum_{j}v^{2}_{j} +{\tilde R}.$$
Denote 
$$\hbox{exp}\left\{ \sum_{jk} \tau t v_{jk}(p)E^{+}_{j}E^{-}_{k}\right\}
\phi(q,p)$$ by $U,$  then $H=\Phi U.$ By the proposition 3 in \S 2, 
Weizenbock formula and a popular formula
$$\Delta_0 (\Phi U)=(\Delta_0 \Phi)U +2\sum_i (\nabla_{E_i}\Phi)
(\nabla_{E_i}U) +\Phi(\Delta_0 U),$$
we have 
$$\begin{array}{rl}
(\frac{\partial}{\partial \tau}+\Box_t)(\Phi U)=&
-[(\Delta_0 -\sum_i \frac{\partial^{2}}{\partial y_{i}^{2}}) \Phi ] U
-2\sum_i (\nabla_{E_i}\Phi)(\nabla_{E_i}U)\\
&+ I_1 +I_2 + I_3 + I_4,
\end{array}$$
where
$$\begin{array}{rl}
I_1&=\Phi(\frac{\partial}{\partial \tau}U -\sum t v_{jk}E_{j}^{+}E_{k}^{-}U)\\
I_2&=t^2[\sum_j v_j^2 - \sum_j (v_{j}(p)+\sum_{k}v_{jk}(p) y_k )^2]\Phi U,\\
I_3&={\tilde R}\Phi U,\\
I_4&=\Phi \Delta_0 U.
\end{array}$$
Note that the lemma 7 still holds if we replace $H$ by $\Phi.$ We write 
the right-hand side of the inequality of the lemma 7 as
$$c_0 Q(c_1 \tau, q,p)\hbox{exp}
\left\{-\frac{\tau t^2}{c_1}v(p)^2 \right\}=
 \sqrt{2}^{n}  c_0 Q(2c_1 \tau, q,p)\hbox{exp}
\left\{-\frac{\tau t^2}{2c_1}v(p)^2 \right\} \xi,$$
where 
$$\xi = \hbox{exp}
\left\{-\frac{1}{8c_1 \tau}\rho^{2}(q,p)-\frac{\tau t^2}{2c_1}v(p)^2 \right\}.
$$   
If we prove
$$|\xi \Phi^{-1}(\frac{\partial}{\partial \tau }+\Box_t )(\Phi U)|<
\hbox{const.}
(\sqrt{\tau }t+1),$$
then
$$\begin{array}{rl}
 |(\frac{\partial}{\partial \tau }+\Box_t )(\Phi U)|&=
|\Phi^{-1}(\frac{\partial}{\partial \tau }+\Box_t )(\Phi U)|\cdot |\Phi| \\
&\leq |\xi \Phi^{-1} (\frac{\partial}{\partial \tau }+\Box_t )(\Phi U)|
{\sqrt{2}}^n c_0 Q(2c_1 \tau,q,p)\hbox{exp}\{-\frac{\tau t^2}{2c_1}v(p)^2 \} \\
&\leq \hbox{const.}c_0 Q(2c_1 \tau,q,p)\hbox{exp}\{-\frac{\tau t^2}{2c_1}v(p)^2 \}
(\sqrt{\tau }t +1).
\end{array}$$
It implies Lemma 8 is correct. So we check 6 terms in the expression of
$(\frac{\partial}{\partial \tau }+\Box_t )(\Phi U)$ 
along this way in order to enable the correctness of Lemma 8. Note that for any 
$m_1, m_2 > 0$
$$\left(\frac{\rho^{2}(q,p)}{\tau }\right)^{m_1}
\left(\tau t^2 v(p)^2
\right)^{m_2}\xi \leq c,$$
where the constant c depends only on $m_1,m_2 .$
So we have
$$\begin{array}{rl}
|\xi \Phi^{-1} I_1|&\leq \xi |\sum t(v_{jk}-v_{jk}(p))E_{j}^{+}E_{k}^{-}| \\
&\leq \hbox{const.}t|y|\xi \\
&=\hbox{const.}t \sqrt{\tau}\left(\frac{|y|}{\sqrt{\tau}}\right)\xi \\
&\leq \hbox{const.} t \sqrt{\tau},
\end{array}$$
$$\begin{array}{rl}
|\xi \Phi^{-1}I_2 |&\leq t^2 |\sum (v_{j}(p)+v_{jk}(p)y_k +O(|y|^2 ))^2 
-\sum(v_{j}(p) +v_{jk}(p)y_{k})^2 |\xi \\
&\leq \hbox{const.} t^2 (|v(p)|\cdot |y|^2 + |y|^3 )\xi \\
&= \hbox{const.}\left[ t \sqrt{\tau}(\sqrt{\tau}t |v(p)|)
(\frac{|y|^2}{\tau}) +
t^2 \sqrt{\tau}^3 (\frac{|y|}{\sqrt{\tau}})^3\right]\xi \\
& \leq \hbox{const.} t \sqrt{\tau},\end{array}$$
where $ ``\hbox{const.}'' $ means a constant, which does not depend on 
$\tau,t,q$ and $p$. Further it is easy to see that
$$|\xi \Phi^{-1} I_3| \leq \hbox{const.},\qquad 
|\xi \Phi^{-1} I_4| \leq \hbox{const.}.$$
So only two terms $[-(\Delta_0 -\sum \frac{\partial^2}{\partial y_{i}^{2}})
\Phi]U$ and   $[-2\sum (\nabla_{E_i}\Phi )(\nabla_{E_i}U)]$
are left alone. From the equalities
$$\begin{array}{rl}
\Phi_i &= \nabla_{E_i}\Phi =H^{ij}\frac{\partial}{\partial y_j}\Phi \\
&=\Phi H^{ij}\left\{ -\frac{1}{2\tau}\left(\frac{\Theta \cosh}{\sinh}
\right)_{jk}y_k -2\tau \left(\frac{\cosh-1}{\Theta \sinh}\right)_{jk}
(Ba^{*})_k \right\}
\end{array}$$
it follows
$$|\xi \Phi^{-1}\Phi_i |\leq \hbox{const.}
\left(\frac{|y|}{\tau} +\tau t^2\right) \xi.$$
Note that $\phi (q,p)=1$ for $q$ , which is near $p$, and 
$\hbox{exp}\left\{\sum \tau t v_{jk}(p)E_{j}^{+}E_{k}^{-}\right\}$ 
does not depend on $y$, so due to $|\Gamma_{ij}^{k}| \leq \hbox{const.} |y|$
we have 
$$|U_i |\leq \hbox{const.}|y|.$$
Therefore
$$|\xi \Phi^{-1}\Phi_i U_i |\leq \hbox{const.}\left(\frac{|y|^2}{\tau}+
\tau t^2 |y|\right)\xi \leq \hbox{const.}(\sqrt{\tau}t +1).$$
By using the equalities
$$\begin{array}{rl}
\sum_i (\Phi_{ii} -\frac{\partial^2}{\partial y_{i}^{2}}\Phi)&
=\Phi (H^{ik}H^{ij} -\delta_{jk})(-\frac{1}{2\tau})
\left(\frac{\Theta \cosh}{\sinh}\right)_{jk} \\
&+\Phi H^{il}\frac{\partial H^{ij}}{\partial y_l}(-\frac{1}{2\tau})
\left(\frac{\Theta \cosh}{\sinh}\right)_{jk}y_k \\
&-\Phi \Gamma_{ii}^{l}H^{lj}\left\{ -\frac{1}{2\tau}
\left(\frac{\Theta \cosh}{\sinh}
\right)_{jk}y_k -2\tau \left(\frac{\cosh-1}{\Theta \sinh}\right)_{jk}
(Ba^{*})_k \right\}
\end{array}$$
and 
$$H^{ij}=\delta_{ij}+\frac{1}{6}\sum R_{ijkl}(p)y_k y_l+\ldots$$
we get
$$|\xi \Phi^{-1}\sum_i (\Phi_{ii} -\frac{\partial^2}{\partial y_{i}^{2}}
\Phi)| \leq \hbox{const.}\left(\frac{|y|^2}{\tau}+\tau t^2 |y|\right)
\xi \leq \hbox{const.}(1+\sqrt{\tau}t).$$
Therefore the lemma is proved.

\section{ Two lemmas}

{\bf Lemma A} \ For a fixed sufficient small $\epsilon >0$, and for any
$c_{1}, c_{2}>0$
with $c_{1}<c_{2}$, there exists a constant $c=c(c_{1},c_{2}, \epsilon)$ such that for 
any $q,p\in M,\nu,\tau >0$ with $ 0<\nu<\tau$, we have
$$\displaystyle\int _{\rho(q,z)< \epsilon} Q(c_{1}(\tau-\nu),q,z)Q(c_{2}\nu,z,p)dz \leq cQ(c_{2}\tau,q,p),$$
where
$$Q(\alpha,q,z)=\frac{1}{\sqrt{4\pi \alpha}^{n}}\hbox{exp}
\left\{-\frac{\rho^{2}(q,z)}{4\alpha}\right\}.$$

{\bf Proof} \ It is equivalent to prove
$$\displaystyle\int _{\rho(q,z)<\epsilon} \sqrt{ \frac{\tau}{(\tau-\nu)\nu} }^n
\hbox{exp}\left\{\frac{1}{4}A(\nu,\tau,q,z,p)\right\} dz \leq 
\hbox{const.},$$
where const. means a constant, which does not depend on $\nu,\tau,q,$ and
$p$, and 
$$A(\nu,\tau,q,z,p)=-\frac{\rho^{2}(q,z)}{c_{1}(\tau-\nu)}
-\frac{\rho^{2}(z,p)}{c_{2}\nu}
+\frac{\rho^{2}(q,p)}{c_{2}\tau}.$$

Now we prove the above inequlity in three separate cases.

(i) If $\rho(q,p)\geq 4\epsilon$ and 
$\frac{\tau-\nu}{\tau}\geq \frac{1}{2}$, then due to $\rho (q,z)<\epsilon,$
$$\rho(z,p)\geq \rho(q,p)-\rho(q,z) \geq 4\epsilon - \epsilon=3\epsilon,$$
$$\frac{\tau-\nu}{\tau}\rho(z,p)-\frac{\nu}{\tau}\rho(q,z) \geq 
\frac{1}{2}\rho(z,p)-\frac{1}{2}\rho(q,z)\geq \frac{3\epsilon}{2}
-\frac{\epsilon}{2}=\epsilon,$$
and thus
$$\begin{array}{rl}
A(\nu,\tau,q,z,p)
&\leq-\frac{\rho^{2}(q,z)}{c_{1}(\tau-\nu)}
-\frac{\rho^{2}(z,p)}{c_{2}\nu}
+\frac{\rho^{2}(q,p)}{c_{2}\tau}\\[3mm]
&\leq-\frac{\rho^{2}(q,z)}{c_{2}(\tau-\nu)}
-\frac{\rho^{2}(z,p)}{c_{2}\nu}
+\frac{(\rho(q,z)+\rho(z,p))^{2}}{c_{2}\tau}\\[3mm]
&=-\frac{\tau}{c_{2}(\tau-\nu)\nu}\left(\frac{\nu}{\tau}\rho(q,z)
-\frac{\tau-\nu}{\tau}\rho(z,p)\right)^{2}\\[3mm]
&\leq -\frac{\tau}{c_{2}(\tau-\nu)\nu}\epsilon^{2}.
\end{array}$$
Therefore
$$\sqrt{\frac{\tau}{(\tau-\nu)\nu}}^{n}\hbox{exp}\frac{A}{4}
\leq\sqrt{\frac{\tau}{(\tau-\nu)\nu}}^{n}\hbox{exp}\{-\frac{\tau}
{4c_{2}(\tau-\nu)\nu}\epsilon^{2}\}\leq
 \hbox{const.}$$

(ii) If $\rho(q,p)\geq 4\epsilon$ and $\frac{\tau-\nu}{\tau}<\frac{1}{2}$,then
$\frac{\nu}{\tau} \geq \frac{1}{2}$ and thus
$$\begin{array}{rl}
\sqrt{\frac{\tau}{(\tau-\nu)\nu}}^{n}\hbox{exp}\frac{A}{4}
&\leq \sqrt{\frac{2}{\tau-\nu}}^{n}\hbox{exp}\left\{
  -(\frac{1}{c_1}-\frac{1}{c_2})\frac{\rho^{2}(q,z)}{4(\tau-\nu)}
  +\frac{1}{4}(-\frac{\rho^{2}(q,z)}{c_{2}(\tau-\nu)}\right. 
\\ &\,\,\,\,\,\,\,\, \left. 
-\frac{\rho^{2}(z,p)}{c_{2}\nu}
+\frac{\rho^{2}(q,p)}{c_{2}\tau})       \right\}\\[3mm]
&\leq \sqrt{\frac{2}{\tau-\nu}}^{n}\hbox{exp}\left\{
  -(\frac{1}{c_1}-\frac{1}{c_2})\frac{\rho^{2}(q,z)}{4(\tau-\nu)}
      -\frac{\tau}{4c_{2}(\tau-\nu)\nu}\left(\frac{\nu}{\tau}\rho(q,z)
-\frac{\tau-\nu}{\tau}\rho(z,p)\right)^{2}           
      \right\}\\[3mm]
&\leq \sqrt{\frac{2}{\tau-\nu}}^{n}\hbox{exp}\left\{
  -(\frac{1}{c_1}-\frac{1}{c_2})\frac{\rho^{2}(q,z)}{4(\tau-\nu)}\right\},
\end{array}$$
and
$$\begin{array}{rl}
\displaystyle\int _{\rho(q,z)<\epsilon} \sqrt{ \frac{\tau}{(\tau-\nu)\nu} }^n
\hbox{exp}\left\{\frac{A}{4}\right\} dz
\leq &
\displaystyle\int _{\rho(q,z)<\epsilon} \sqrt{ \frac{2}{(\tau-\nu)} }^n
\hbox{exp}\left\{-\frac{\rho^2(q,z)}{4c_0 (\tau-\nu)}\right\} dz\\
\leq & \hbox{const.},
\end{array}$$
where $c_0 =(\frac{1}{c_1}-\frac{1}{c_2})^{-1}.$
Therefore the lemma is true in this case.

(iii) Now we consider the case when  $\rho(q,p)<4\epsilon$. Let
$$\lambda=\frac{\nu}{\tau},\qquad \mu=\frac{\tau - \nu}{\tau}.$$
By using a reasonning in the proof of case (i),
$$\begin{array}{rl}
\lambda \rho^2 (q,z)+\mu\rho^2 (z,p) -\lambda\mu\rho^2 (q,p)&\geq 
\lambda \rho^2 (q,z)+\mu\rho^2 (z,p) -\lambda\mu(\rho (q,z)+\rho(z,p))^2 \\
&=(\lambda \rho(q,z)-\mu\rho(z,p))^2 \\
&\geq 0,
\end{array}$$
 we can let
$$W=\sqrt{\lambda \rho^2 (q,z)+\mu\rho^2 (z,p) -\lambda\mu\rho^2 (q,p)}.$$
And let $o$ be a point on the geodesic joining $p$ and $q$ such that
$$\frac{\rho (o,p)}{\rho (q,p)}=\lambda.$$
Without of loss of generalities, we assume $n=2$.
By the last
 theorem in the appendix we know that for a sufficeinte small
$ \epsilon>0,$ 
$$W^2 \geq \frac{1}{2}\rho^2 (o,z).$$
Then
$$\begin{array}{rl}
 \sqrt{\frac{\tau}{(\tau-\nu)\nu}}^{n} \hbox{exp}\frac{A}{4}&
\leq \sqrt{\frac{\tau}{(\tau-\nu)\nu}}^{n} \hbox{exp}\{
-\frac{\tau}{4c_2 (\tau-\nu)\nu}W^2\} \\
&\leq \sqrt{\frac{\tau}{(\tau-\nu)\nu}}^{n} \hbox{exp}\{
-\frac{\tau}{8c_2 (\tau-\nu)\nu}\rho^2 (o,z)\}.
\end{array}$$
Choose a geodesic coordinate system centering at $o$, the coordinates
of a point is $(\rho,\theta),$  where $\theta \in S^{1}.$
Then the volume measure $dz$ satisfies
$$dz\leq \hbox{const.}|\rho d\rho d\theta|.$$
 And then
$$\sqrt{\frac{\tau}{(\tau-\nu)\nu}} \hbox{exp}\frac{A}{4}dz \leq\hbox{const.}
  \sqrt{\frac{\tau}{(\tau-\nu)\nu}} \hbox{exp}\{
-\frac{1}{8c_2}\frac{\tau}{(\tau-\nu)\nu}\rho^2 \}\rho d\rho d\theta.$$
Therefore we finish the proof of Lemma A.

{\bf Lemma B}\ Let $\epsilon$ be small enough. For any $s_{0}>0, c>0$ there
exists a constant $h(s_0,c)$ such that for any ${\tilde c}\geq h(s_0,c)$ and
$t,\tau_1,\tau_2>0$ with $\tau_{1}t, \tau_{2}t <s_0$, and $q,p\in M$ with
$\rho(q,p)<\epsilon$, we have
$$\hbox{exp}\left\{
-\frac{\rho^{2}}{4c\tau_{1}}-\frac{\tau_{1}t^{2}}{c}v(p)^{2}\right\}
\hbox{exp}\left\{  -\frac{\tau_{2}t^{2}}{\tilde c}v(p)^{2}\right\}
\leq \hbox{exp}\left\{-\frac{\tau_{2}t^{2}}{\tilde c}v(q)^{2}\right\}$$
and
$$\hbox{exp}\left\{
-\frac{\rho^{2}}{4c\tau_{1}}-\frac{\tau_{1}t^{2}}{c}v(p)^{2}\right\}
\hbox{exp}\left\{-\frac{\tau_{2}t^{2}}{\tilde c}v(q)^{2}\right\}
\leq \hbox{exp}\left\{-\frac{\tau_{2}t^{2}}{\tilde c}v(p)^{2}\right\},$$
where $\rho=\rho(q,p)$.

{\bf Proof} \ It is equivalent to prove
$$\hbox{exp}\left\{
-\frac{\rho^{2}}{4c\tau_{1}}-\frac{\tau_{1}t^{2}}{c}v(p)^{2}\right\}
\hbox{exp}\left\{\frac{\tau_{2}t^{2}}{\tilde c}|v(q)^{2}-v(p)^{2}|\right\}
\leq 1,$$
or 
$$
-\frac{\rho^{2}}{4c\tau_{1}}-\frac{\tau_{1}t^{2}}{c}v(p)^{2}
+\frac{\tau_{2}t^{2}}{\tilde c}|v(q)^{2}-v(p)^{2}|\leq 0.$$
From Taylor expansion of $v(q)^{2}$
$$v(q)^{2}=v(p)^{2}+2v(p)v^{'}(p)\rho +\ldots,$$
it follows that there exists a constant $k$, which does not depend on 
$p$ and $q$, such that 
$$|v(q)^{2}-v(p)^{2}|\leq k(|v(p)|\rho + \rho^{2}).$$
So
$$\begin{array}{rl}
-\frac{\rho^{2}}{4c\tau_{1}}&-\frac{\tau_{1}t^{2}}{c}v(p)^{2}
+\frac{\tau_{2}t^{2}}{\tilde c}|v(q)^{2}-v(p)^{2}|\leq 
-\frac{\rho^{2}}{4c\tau_{1}}-\frac{\tau_{1}t^{2}}{c}v(p)^{2}
+\frac{k\tau_{2}t^{2}}{\tilde c} (|v(p)|\rho +\rho^{2})   \\[3mm]
&\leq -\frac{\rho^{2}}{4c\tau_{1}}-\frac{\tau_{1}t^{2}}{c}v(p)^{2}
+\frac{k\tau_{2}t^{2}}{\tilde c}( \frac{2c}{t}\sqrt{\frac{v(p)^{2}\tau_{1}t^{2}}
{c}}\sqrt{\frac{\rho^{2}}{4c\tau_{1}}} )
+\frac{4kc\tau_{1}\tau_{2}t^{2}}{\tilde c}\frac{\rho^{2}}{4c\tau_{1}}\\[3mm]
&\leq -\frac{\rho^{2}}{4c\tau_{1}}-\frac{\tau_{1}t^{2}}{c}v(p)^{2}
  +\frac{k\tau_{2}tc}{\tilde c}
(\frac{\rho^{2}}{4c\tau_{1}}+\frac{\tau_{1}t^{2}}{c}v(p)^{2})
+\frac{4kc\tau_{1}\tau_{2}t^{2}}{\tilde c} \frac{\rho^{2}}{4c\tau_{1}}\\[3mm]
&\leq -\frac{\rho^{2}}{4c\tau_{1}}(1-\frac{ks_{0}c}{\tilde c}-\frac{4ks_{0}^{2}c}
{\tilde c})
-\frac{\tau_{1}t^{2}}{c}v(p)^{2}(1-\frac{ks_{0}c}{\tilde c}).
\end{array}$$
Therefore if  choose $h(s_{0},c)\geq (ks_{0} +4ks_{0}^{2})c$, then we have
$$
-\frac{\rho^{2}}{4c\tau_{1}}-\frac{\tau_{1}t^{2}}{c}v(p)^{2}
+\frac{\tau_{2}t^{2}}{\tilde c}|v(q)^{2}-v(p)^{2}|\leq 0.$$
So lemma B is true.

\section {Levi iteration}
{\bf Definition 9}\ Suppose we are given $H(\tau,q,p,t).$ By the following
procedure we construct   $K_m (\tau,q,p,t), $   $m>0,$ and 
$K (\tau,q,p,t).$
$$\begin{array}{rl}
K_0 (\tau,q,p,t)&=(\frac{\partial}{\partial \tau} +\Box_t )H(\tau,q,p,t), \\
K_{m+1} (\tau,q,p,t)&=\displaystyle\int_{0}^{\tau}d\nu \displaystyle\int
K_0 (\tau-\nu,q,z,t)K_m (\nu,z,p,t)dz,\quad \forall m\geq 0, \\
K(\tau,q,p,t)&=\displaystyle\sum_{m=0}^{\infty}(-1)^{m+1}
K_m (\tau,q,p,t).
\end{array}$$
The above procedure is called Levi iteration.

Of course, the first question for Levi iteration is wheather the series,
which defines $K(\tau,q,p,t)$, converges. We will show it does if
$H(\tau,q,p,t)$ is defined by Definition 6.

{\bf Lemma 10}\ Choose $\phi(q,p)$ properly such that the supports of 
$H(\tau,q,p,t)$ and $(\frac{\partial}{\partial \tau} +\Box_t )H(\tau,q,p,t)$
are contained in
$\{(\tau,q,p,t)|\rho(q,p)<\epsilon\}.$ Then for sufficient small $\epsilon$ 
and fixed $\tau_0 ,s_0 >0$ there exist $c_0,c_1$ such that Lemma 7, Lemma 8 and the following inequalities hold.
$$\begin{array}{rl}
&|H(\tau,q,p,t)|\leq c_0 Q(c_1 \tau,q,p)\hbox{exp}\left\{-\frac{\tau t^2}{c_1}
v(q)^2 \right\},\\
&|(\frac{\partial}{\partial \tau} +\Box_t )H(\tau,q,p,t)|
\leq c_0 (\sqrt{\tau}t +1)Q(c_1 \tau,q,p)
\hbox{exp}\left\{-\frac{\tau t^2}{c_1}v(q)^2\right\}.
\end{array}$$

{\bf Proof}\ The lemma is trivial due to Lemma7, Lemma 8 and Lemma B.

{\bf Lemma 11}\ Choose $\epsilon,\tau_0,s_0$ as in Lemma 10 and let $c_0,c_1$
be given in Lemma 10. Then for ${\tilde c_1}>\hbox{Max}\{h(s_0,2c_1), 2c_1 
,\epsilon\},$
where $h(s_0, 2c_1,\epsilon)$ is given in Lemma B, the following inequalities hold
$$
|K_m (\tau,q,p,t)|\leq \left(\frac{\tilde c_1}{c_1}\right)^{\frac{n}{2}}
c_0 (\sqrt{2}^{n} c_0 c)^m (\sqrt{\tau}t +1)^{m+1}
\frac{\tau^m }{m!}Q({\tilde c_1}\tau,q,p)\hbox{exp}\left\{
-\frac{\tau t^2}{\tilde c_1}v(q)^2\right\},$$
 where $ m\geq 0,$ 
and  $c=c(2c_1,{\tilde c_1}),$ which is given by Lemma A.

{\bf Proof} Let
$$ a_{0}=\left(\frac{\tilde c_1}{c_1}\right)^{\frac{n}{2}} c_{0},$$
then by Lemma 10
$$\begin{array}{rcl}
|K_{0}(\tau,q,p,t)|&=&
|(\frac{\partial}{\partial \tau} +\Box_t )H(\tau,q,p,t)| \\
&\leq& c_{0}(\sqrt{\tau} t +1) \frac{1}{{\sqrt{4 \pi c_1 \tau}}^n }
\hbox{exp}\left\{-\frac{\rho^{2}(q,p)}{4c_1 \tau}
-\frac{\tau t^2}{ c_1}V(q) \right\} \\
&\leq&\left(\frac{\tilde c_1}{c_1}\right)^{\frac{n}{2}} c_{0}
   (\sqrt{\tau} t +1) \frac{1}{{\sqrt{4 \pi {\tilde c_1} \tau}}^n }
\hbox{exp}\left\{-\frac{\rho^{2}(q,p)}{4{\tilde c_1} \tau}
-\frac{\tau t^2}{\tilde c_1} V(q) \right\} \\
&=&
a_0 (\sqrt{\tau}t +1) Q({\tilde c_1}\tau,q,p)\hbox{exp}\left\{
-\frac{\tau t^2}{\tilde c_1}v(q)^2 \right\}.
\end{array}$$
Define
$$a_{m+1}=(\sqrt{2}^{n}c_0 c)a_{m}=
\left(\frac{\tilde c_1}{c_1}\right)^{\frac{n}{2}} c_0
(\sqrt{2}^{n}c_0 c)^{m+1}.$$
We are going to check the following equalities
$$
|K_m (\tau,q,p,t)|\leq a_m  (\sqrt{\tau}t +1)^{m+1}
\frac{\tau^m }{m!}Q({\tilde c_1}\tau,q,p)\hbox{exp}\left\{
-\frac{\tau t^2}{\tilde c_1}v(q)^2\right\},$$
by induction on $m$. Suppose the inequality is true for $m$,
by using an inequality in Lemma 10 and a fact
$$\hbox{Support}(K_0 (\tau,q,p,t))\subset
\{(\tau,q,p,t)|\rho(q,p)<\epsilon\}$$
we have
$$\begin{array}{rcl}
|K_{m+1}|&=&|\displaystyle\int_{0}^{\tau}\displaystyle\int_{M}
K_0 (\tau-\nu,q,z,t)K_m (\nu,z,p,t)dzd\nu|\\
&=&|\displaystyle\int_{0}^{\tau}\displaystyle\int_{\rho (z,q) < \epsilon}
K_0 (\tau-\nu,q,z,t)K_m (\nu,z,p,t)dzd\nu|\\
&\leq&\displaystyle\int_{0}^{\tau}\displaystyle\int_{\rho (z,q) < \epsilon}
c_0 a_{m}(\sqrt{\tau - \nu}t +1) (\sqrt{\nu}t +1)^{m + 1}
\frac{\nu^m}{m!}{\tilde Q}dz d\nu\\
&\leq&\displaystyle\int_{0}^{\tau}\displaystyle\int_{\rho (z,q) < \epsilon}
c_0 a_{m} (\sqrt{\tau}t +1)^{m + 2}
\frac{\nu^m}{m!}{\tilde Q}dz d\nu\\,
\end{array}$$
where
$${\tilde Q}=
 Q(c_1 (\tau-\nu),q,z)\hbox{exp}\left\{-\frac{(\tau-\nu) t^2}{c_1}
 v(z)^2 \right\}
 Q({\tilde c_1} \nu,z,p)\hbox{exp}\left\{-\frac{\nu t^2}{\tilde c_1}
 v(z)^2 \right\}.$$
 By Lemma B we have
$$\begin{array}{rcl}
{\tilde Q}&=&
\sqrt{2}^{n}Q(2c_1 (\tau-\nu),q,z)Q({\tilde c_1}\nu,z,p)
\hbox{exp}\left\{-\frac{\rho^2(q,z)}{8(\tau-\nu)c_1}
-\frac{(\tau-\nu)t^2}{2c_1}v(z)^2\right\}\\   
&&\qquad \hbox{exp}\left\{-\frac{(\tau-\nu)t^2}{2c_1}v(z)^2\right.  \\ 
&&\,\,\,\,\,\,\,\, \left.
  -\frac{\nu t^2}{\tilde c_1}v(z)^2 \right\}\\
&=&\sqrt{2}^{n}Q(2c_1 (\tau-\nu),q,z)Q({\tilde c_1}\nu,z,p)
\hbox{exp}\left\{-\frac{\rho^2(q,z)}{8(\tau-\nu)c_1}
-\frac{(\tau-\nu)t^2}{2c_1}v(z)^2\right\}\\
&&\qquad \hbox{exp}\left\{-\frac{(\tau-\nu)t^2}{2{ c_1}}v(z)^2
  -\frac{\nu t^2}{\tilde c_1}v(z)^2 \right\}\\
&\leq &\sqrt{2}^{n}Q(2c_1 (\tau-\nu),q,z)Q({\tilde c_1}\nu,z,p)
\hbox{exp}\left\{-\frac{\rho^2(q,z)}{8(\tau-\nu)c_1}
-\frac{(\tau-\nu)t^2}{2c_1}v(z)^2\right\}\\
&&\,\,\,\,\,\,\,\,
\hbox{exp}\left\{-\frac{\tau t^2}{\tilde c_1}v(z)^2  \right\}\\
&\leq& \sqrt{2}^{n}Q(2c_1 (\tau-\nu),q,z)Q({\tilde c_1}\nu,z,p)
\hbox{exp}\left\{-\frac{\tau t^2}{\tilde c_1}v(q)^2  \right\}.
\end{array}$$
And by Lemma A we also have
$$\displaystyle\int_{\rho (z,q) < \epsilon}{\tilde Q}dz
\leq {\sqrt 2}^n c 
Q({\tilde c_1}\tau,q,p)
\hbox{exp}\left\{-\frac{\tau t^2}{\tilde c_1}v(q)^2  \right\}.$$
Therefore,
$$\begin{array}{rcl}
|K_{m+1}|&\leq&\displaystyle \int_{0}^{\tau} (c_0 c {\sqrt 2}^n) a_m
 ({\sqrt \tau}t +1)^{m+2}\frac{\nu^{m}}{m!}d\nu \cdot
Q({\tilde c_1}\tau,q,p)
\hbox{exp}\left\{-\frac{\tau t^2}{\tilde c_1}v(q)^2  \right\} \\
&=&  a_{m+1}
 ({\sqrt \tau}t +1)^{m+2}\frac{\tau^{m+1}}{(m+1)!}
Q({\tilde c_1}\tau,q,p)
\hbox{exp}\left\{-\frac{\tau t^2}{\tilde c_1}v(q)^2  \right\}. 
\end{array}$$
The lemma is proved.

{\bf Lemma 12}\quad There holds
$$\begin{array}{rcl}
\sum_{m=0}^{\infty}|K_m (\tau,q,p,t)|&\leq&
\left(\frac{\tilde c_1}{c_1}\right)^{\frac{n}{2}} c_0 (\sqrt{\tau}t+1)
\hbox{exp}\left\{\sqrt{2}^n c_0 c(\sqrt{\tau_0} s_0
+\tau_0)\right\}\\
&&
Q({\tilde c_1}\tau,q,p)
\hbox{exp}\left\{-\frac{\tau t^2}{\tilde c_1}v(q)^2\right\}.
\end{array}$$

{\bf Proof}\quad It is a trivial corollary of Lemma 11.

{\bf Proposition 13}\quad There holds
$$G(\tau,q,p,t)= H(\tau,q,p,t)
 +\displaystyle\int_{0}^{\tau}d\nu \displaystyle\int_{M}
H(\tau-\nu,q,z,t)K(\nu,z,p,t)dz .$$

{\bf Proof}\quad $K(\tau,q,p,t)$ is well-defined due to Lemma 12. The 
right-hand side of the above equality is also well-defined by using Lemma 12
 and Lemma A. By a necessary routine check the right-hand side is indeed
a fundamental solution of $\frac{\partial}{\partial \tau}+\Box_t.$
So the proposition is true.

\section{Proof of Theorem 1}

Due to 
$$\Theta^\# =\sqrt{4\tau^2 B^* B}=\sqrt{4\tau^2 t^2 A(p)A(p)^*}=\theta$$
we have 
$H(\tau,p,p,t)=\phi_0 (\tau,t,p).$ Then by Proposition 13, Theorem 1 is 
equivalent to 
$$(s-\lim) \displaystyle\int_{M}dp\displaystyle\int_{0}^{\tau}d\nu
\displaystyle\int_{M}H(\tau-\nu,p,z,t)K(\nu,z,p,t)dz=0.$$

By lemma 10, lemma 12, and lemma B, which is used in the same way 
as in proving Lemma 11, we can get
$$\begin{array}{rl}
(i) \quad |H(\tau-\nu,p,z,t)K(\nu,z,p,t)|\leq \hbox{const.}&
(\sqrt{\tau}t+1)Q(2c_1 (\tau-\nu),p,z)Q({\tilde c_1}\nu,z,p) \\
&\hbox{exp}\left\{-\frac{\tau t^2}{\tilde c_1}v(p)^2 \right\},\\
(ii) \quad |H(\tau-\nu,p,z,t)K(\nu,z,p,t)|\leq \hbox{const.}&
(\sqrt{\tau}t+1)Q(2c_1 (\tau-\nu),p,z)Q({\tilde c_1}\nu,z,p) \\
&\hbox{exp}\left\{-\frac{\tau t^2}{\tilde c_1}v(z)^2 \right\},
\end{array}$$
 where the $\hbox{const.}$ is a constant, which does not depend
 on $\tau, t, q, p, z.$

Let $N_{\epsilon}$ be an $\epsilon$-neighbourhood of the zero set of $V,$ 
and let $$\delta=\hbox{Min}\{|V(p)|; p\not\in N_{\epsilon}\}.$$
Thus if $z$ is not in $N_{\epsilon}$, then by (ii) we have
$$\begin{array}{rl}
 |H(\tau-\nu,p,z,t)K(\nu,z,p,t)|\leq \hbox{const.}&
(\sqrt{\tau}t+1)Q(2c_1 (\tau-\nu),p,z) \\
&Q({\tilde c_1}\nu,z,p)
\hbox{exp}\left\{-\frac{\tau t^2 \delta^2}{\tilde c_1} \right\},
\end{array}$$
$$\begin{array}{rl}
|\displaystyle\int_{M} dp\displaystyle\int_{0}^{\tau}d\nu&
\displaystyle\int_{M-N_{\epsilon}}H(\tau-\nu,p,z,t)K(\nu,z,p,t)dz|
\leq      \\
&\leq\hbox{const.}
\hbox{exp}\left\{-\frac{\tau t^2 \delta^2}{\tilde c_1} \right\}
\displaystyle\int_{0}^{\tau}(\sqrt{\tau}t+1)d\nu  \\
&\qquad
{\displaystyle\int\!\!\displaystyle\int}_{\rho(p,z)<\epsilon}
Q(2c_1 (\tau-\nu),p,z)Q({\tilde c_1}\nu,z,p)dpdz \\
&\leq \hbox{const.}
\hbox{exp}\left\{-\frac{s\delta^2}{\tilde c_1}t\right\}
(\sqrt{\tau}s+\tau)\displaystyle\int_{M}Q({\tilde c_1}\tau,p,p)dp \\
&\leq \hbox{const.}\frac{1}{\sqrt{\tau}^n}
\hbox{exp}\left\{-\frac{s^2 \delta^2}{\tilde c_1}\frac{1}{\tau}\right\},
\end{array}$$
which implies
$$(s-\lim) |\displaystyle\int_{M}dp\displaystyle\int_{0}^{\tau}d\nu
\displaystyle\int_{M-N_{\epsilon}}H(\tau-\nu,p,z,t)K(\nu,z,p,t)dz|=0.$$
Similarly, for $p\not\in N_{\epsilon},$ we also have 
$$(s-\lim) |\displaystyle\int_{M-N_{\nu}}dp\displaystyle\int_{0}^{\tau}d\nu
\displaystyle\int_{M}H(\tau-\nu,p,z,t)K(\nu,z,p,t)dz|=0.$$
Suppose $0\in \hbox{Zero}(V),$ 
choose a normal coordinate system centering at $0$, and an orthonormal
moving frame as in \S 3. Let the coordinates 
of $p$ and $z$ are $(x_1, \cdots,x_n), (z_1,\cdots,z_n) $ respectively.
 And let 
$Y=(y_1,\cdots,y_n)$ be defined by
$$\hbox{exp}_{z}(\sum y_i E_{i}(z))=p.$$ Then 
$$\begin{array}{rl}
H(\tau,p,z,t)=&\Phi (\tau,Y,t(v_1 (z),\cdots,v_n (z)),tB(z)) \\
&\hbox{exp}\left\{\displaystyle\sum \tau t v_{jk}(z)E_{j}^{+}E_{k}^{-}
\right\}\phi(p,z),\end{array}$$
and
$$|H(\tau,p,z,t)|\leq \hbox{const.}\Phi(\tau,Y,t(v_1(z),\cdots,v_n(z)),
tB(z)).$$
Note that
$$Y=(X-Z)+\ldots,$$
$$\begin{array}{rl}
(v_1(z),\cdots,v_n(z))=&(z_1,\cdots,z_n)\left(
   \begin{array}{ccc}
     v_{11}(0)&\cdots&v_{n1}(0) \\
     \vdots&&\vdots \\
     v_{1n}(0)&\cdots&v_{nn}(0)
     \end{array}   \right)+\ldots \\
=&ZB(0) + \ldots
\end{array}$$
and
$$\begin{array}{rl}
\Phi_0 (\tau,X-Z,tZB(0),tB(0))=&\frac{1}{\sqrt{4 \pi \tau}^{n}}
  \sqrt {\det\left(\frac{\Theta}{\sinh}
                   \right)
         }
 \exp\left\{     -\frac{1}{8\tau} (X+Z)
                      \frac{\Theta(\cosh-1)}{\sinh}
                      (X+Z)^{\ast}\right.             \\[3mm]
 &\quad-\frac{1}{8\tau} (X-Z)
                     \frac{\Theta(\cosh+1)}{\sinh}
                     (X-Z)^{\ast}
         \left. \right\}\\
\end{array}$$

Now we estimate  $\Phi(\tau,Y,t(v_1(z),\cdots,
v_n(z)),tB(z))$ and $\Phi(\tau,X-Z,tZB(0),tB(0)).$

{\bf Lemma 14}\quad For fixed $s>0,$ there exist $\epsilon, \eta >0$
such that if $\rho(p,z)\leq \epsilon,$then
$$\begin{array}{rl}
(i)\quad \Phi(\tau,X-Z,tZB(0),tB(0))
\leq \hbox{const.}\frac{1}{\sqrt{4\pi\tau}^n}
\exp\left\{ \right. & -\frac{1}{4\tau\eta}[(X+Z)(X+Z)^*  \\
&+(X-Z)(X-Z)^*]\left. \right\}.\\
(ii)\quad \Phi(\tau,Y,t(v_1(z), \cdots,v_n(z)),tB(z))
\leq \hbox{const.}\frac{1}{\sqrt{4\pi\tau}^n}
&\exp\left\{-\frac{1}{8\tau\eta}[(X+Z)(X+Z)^*\right. \\
+(X-Z)(X-Z)^*]\left. \right\},
\end{array}$$
where the $\hbox{const.}$ depends on $s$.

{\bf Proof}  Denote
$$\begin{array}{rl}
A_0 = &-\frac{1}{4\tau} 
       (X-Z) \frac{\Theta_0 \cosh_0}{\sinh_0}(X-Z)^{\ast} \\
  &-2\tau t^2 (X-Z) \frac{\cosh_0 -1}{\Theta_0\sinh_0} B(0)
(ZB(0))^{\ast}\\
  &-2 \tau  t^2 ZB(0)\frac{\cosh_0^\# -1}{\Theta_0^\#\sinh_0^\#} 
(ZB(0))^{\ast},\end{array}$$
and
$$\begin{array}{rl}
A=
  &-\frac{1}{4\tau} Y \frac{\Theta\cosh}{\sinh}Y^{\ast}\\
&  -2\tau t^2 Y \frac{\cosh -1}{\Theta\sinh} B(z)
v(z)^{\ast}\\
&  -2 \tau  t^2 v(z)\frac{\cosh^\# -1}{\Theta^\#\sinh^\#} 
v(z)^{\ast},\end{array}$$
where
$$\begin{array}{rl}
v(z)&=(v_1(z),\cdots,v_n(z)),\\
\Theta&=\sqrt{4\tau^2 t^2 B(z)B(z)^{\ast}},\\
\Theta_0 &=\sqrt{4\tau^2 t^2 B(0)B(0)^{\ast}}\\
\cosh_0&=\cosh\Theta_0.
\end{array}$$

It is trivial that there is a positive function $\alpha (\cdot)$ with
$\lim_{\epsilon \to 0}\alpha(\epsilon)=0,$ such that
$$\begin{array}{rl}
|A-A_0|\leq& \hbox{const.}\alpha(\epsilon)\{\frac{1}{4\tau}|X-Z|^2
+2\tau t^2 |X-Z|\cdot|Z| +2\tau t^2 |Z|^2\}\\
\leq& \hbox{const.}\alpha(\epsilon)\{\frac{1}{\tau}|X-Z|^2
+\frac{1}{\tau}|X-Z|\cdot |Z| +\frac{1}{\tau}|Z|^2 \}.\end{array}$$
From 
$$\begin{array}{rl}
&|Z|\leq \frac{1}{2}(|X-Z|+|X+Z|),\\
&|Z|^2\leq \frac{1}{4}(|X-Z|+|X+Z|)^2 \leq \frac{1}{2}
(|X-Z|^2 +|X+Z|^2),\end{array}$$
it follows that
$$|A-A_0|\leq \hbox{const.}\alpha(\epsilon)\{\frac{1}{\tau}|X-Z|^2
+|X+Z|^2\}.$$
By using Proposition 5 we can choose $\eta>0$  such that the following 
inequalities hold.
$$\begin{array}{rl}
A_0 =&-\frac{1}{8}(X-Z)\frac{\Theta_0 (\cosh_0+1)}{\sinh_0}(X-Z)\\
&\quad \quad -\frac{1}{8}(X+Z)\frac{\Theta_0 (\cosh_0-1)}{\sinh_0}(X+Z)\\
\leq& -\frac{1}{4\tau\eta}|X-Z|^2 -\frac{1}{4\tau\eta}|X+Z|^2.
\end{array}$$
So from
$$\frac{\Theta}{\sinh}\geq \hbox{const.}>0$$
we get (i). 
And it is easy to see that
$$\begin{array}{rl}
|\Phi|&\leq 
\hbox{const.}\frac{1}{\sqrt{4\pi \tau}^n}\exp A \leq 
\hbox{const.}\exp\{A-A_0\}\frac{1}{\sqrt{4\pi \tau}^n}\exp A_0 \\
&\leq \hbox{const.exp}\left\{\hbox{const.}\frac{\alpha(\epsilon)}{\tau}
(|X-Z|^2 +|X+Z|^2) -\frac{1}{8\tau\eta}(|X-Z|^2 +|X+Z|^2)\right\} \\
&\qquad \frac{1}{\sqrt{4\pi\tau}^n}\exp
\left\{-\frac{1}{8\tau\eta}(|X-Z|^2 +|X+Z|^2)\right\}.\end{array}$$
Choose $\epsilon$ small enough such that
the term 
$$\left\{\hbox{const.}\frac{\alpha(\epsilon)}{\tau}
(|X-Z|^2 +|X+Z|^2) -\frac{1}{2\tau\eta}(|X-Z|^2 +|X+Z|^2)\right\}$$
in the above inequalities is negative, thus the lemma is true.

Now let us continue to prove the theorem 1.  By using Lemma 12 we have
$$\begin{array}{rl}
|H(\tau-\nu,p,z,t)&K(\nu,z,p,t)|\leq \hbox{const.}\Phi(\tau-\nu,p,z,t)
(\sqrt{\nu}t+1)Q({\tilde c_1}\nu,q,p) \\
\leq&  \hbox{const.}(\sqrt{\tau}t+1)
\frac{1}{\sqrt{\tau-\nu}^n}\frac{1}{\sqrt{\nu}^n} 
\exp\left\{
-\frac{1}{8(\tau-\nu)\eta}(X+Z)(X+Z)^*\right\} \\
& \exp\left\{-\frac{1}{8(\tau-\nu)\eta}(X-Z)(X-Z)^*
\right\} \cdot
\exp \left\{-\frac{1}{4\nu{\tilde c_1}}(X-Z)(X-Z)^*
\right\}.\end{array}$$
Let $W_1=X+Z, \quad W_2=X-Z,$ 
note that
$$\begin{array}{rl}
b_0 \equiv \displaystyle\int
\frac{1}{\sqrt{\tau-\nu}^n}&
\exp\left\{-\frac{1}{8(\tau-\nu)\eta}W_2 W_2^*\right\}
\frac{1}{\sqrt{\nu}^n}
\exp\left\{-\frac{1}{4\nu{\tilde c_1}}W_2 W_2^*\right\}dW_2\\
&\leq \displaystyle\int \frac{1}{\sqrt{(\tau-\nu)\nu}^n}
\exp\left\{-\frac{\tau}{4(\tau-\nu)\nu{\tilde c_2}}
W_2 W_2^*\right\}dW_2\\
&\leq \hbox{const.}\frac{1}{\sqrt{\tau}^n},\end{array}$$
where $${\tilde c_2}=\hbox{Max}\{2\eta,{\tilde c_1}\}.$$
Then 
$$\begin{array}{rl}
|\displaystyle\int_{0}^{\tau}d\nu 
{\displaystyle\int\!\!\displaystyle\int}_{\{\hbox{near\quad} 0\}}&
H(\tau-\nu,p,z,t)K(\nu,z,p,t)dp dz|\\
\leq& \hbox{const.}|\displaystyle\int_{0}^{\tau}d\nu
{\displaystyle\int \!\!\displaystyle\int}_{\{\hbox{near } 0\}}H \cdot K dW_1 
dW_2|\\
\leq&\hbox{const.}\displaystyle\int_{0}^{\tau}(\sqrt{\tau}t+1)d\nu 
\displaystyle\int  b_0
\exp\left\{-\frac{W_1^2}{8(\tau-\nu)\eta}\right\}dW_1 \\
\leq&\hbox{const.}\displaystyle\int_{0}^{\tau}(\sqrt{\tau}t+1)d\nu 
\displaystyle\int\frac{1}{\sqrt{\tau}^n}
\exp\left\{-\frac{W_1^2}{8(\tau-\nu)\eta}\right\}dW_1 \\
\leq&\hbox{const.}\displaystyle\int_{0}^{\tau}(\sqrt{\tau}t+1)d\nu 
\displaystyle\int\frac{1}{\sqrt{\tau}^n}
\exp\left\{-\frac{W_1^2}{4\tau \eta}\right\}dW_1 \\
\leq&\hbox{const.}(\sqrt{\tau}t +1)\tau\\
\leq&\hbox{const.}(\sqrt{\tau}s +\tau) \stackrel{\hbox{s}-\lim}
{\longrightarrow} 0,
\end{array}$$
where  $\{\hbox{near\quad} 0\}$ means that both $p$ and $z$ are near to the 
point $0$.

Suming up the above discussions we get
$$(s-\lim) \displaystyle\int_{M}dp \displaystyle\int_{0}^{\tau}d\nu
\displaystyle\int_{M}
H(\tau-\nu,p,z,t)K(\nu,z,p,t)dz=0.$$

\section{Hopf theorem}

{\bf Hopf Theorem} \ Given a Riemannian manifold $M$ of dim $n$, and a vector 
field $V$ without degenerate zeros, there holds 
$$\chi(M) = \displaystyle\sum_{p \in \hbox{Zero} (V)} \displaystyle\frac {\det(v_{ij}(p))} 
{|\det(v_{ij}(p))|},$$
where $\chi(M)$ is the Euler number of $M$.

{\bf Proof}\quad For the superstructure of $\Lambda^*_p(M)$
$$\Lambda^*_p(M)=\Lambda^{even}_p(M)+\Lambda^{odd}_p(M),$$
define $\hbox{str}$ of a linear map
$$L:\Lambda^*_p(M) \rightarrow \Lambda^*_p(M)$$
by
$$\hbox{str}(L)=\hbox{tr}(L|\Lambda^{even}_p(M))
-\hbox{tr}(L|\Lambda^{odd}_p(M)).$$
Due to the following inequalites
$$|a_1 +\cdots +a_n |\leq |a_1|+\cdots + |a_n|\leq
\sqrt{n}\sqrt{a_1^2 +\cdots +a_n^2},$$
we have
$$|\hbox{str}(G(\tau,p,p,t)-\phi_0 (\tau,t,p))|\leq
\sqrt{n} |G(\tau,p,p,t)-\phi_0 (\tau,t,p)|.$$
Thus the theorem 1 implies
$$\chi (M)=(s-\lim) \displaystyle\int_M \hbox{str}G(\tau,p,p,t)dp=
(s-\lim) \displaystyle\int_M \hbox{str}\phi_0(\tau,t,p)dp.$$
Let $N_{\epsilon}$ be the $\epsilon$-neighbourhood of 
Zero$(V)$ . It is easy to see that there exists an $\delta>0$
such that 
$$\delta \leq \hbox{Min}\{|V(p)|;p\in M-N_{\epsilon}\},$$
and
$$\delta \leq \frac{\cosh\theta -1}{\theta \sinh\theta.}
\hbox{\quad\quad}\forall s\leq s_0$$
Then
$$|\displaystyle\int_{M-N_{\epsilon}} \hbox{str}\phi_0(\tau,t,p)dp|\leq
\hbox{const.}\displaystyle\int_{M-N_{\epsilon}} \frac{1}{\sqrt{4\pi\tau}^n}
e^{-2\tau s \delta^3}dp \stackrel{\hbox{s}-\lim}{\longrightarrow} 0.$$
So
$$(s-\lim) \displaystyle\int_M \hbox{str}\phi_0(\tau,t,p)dp=
(s-\lim)\displaystyle\int_{N_{\epsilon}} \hbox{str}\phi_0(\tau,t,p)dp.$$

Of course, Zero$(V)$ is a finite set $\{p_1,\cdots,p_m\}.$ For each zero 
point $p_{\alpha}$, in its $\epsilon-$neighbourhood $N_{\epsilon}(p_{\alpha})$,
choose a normal coordinate system centering at $p_{\alpha},$ and orthonormal
frames $\{E_1,\cdots,E_n\}$ as before. Denote the normal coordinates of 
$q\in N_{\epsilon}(p_{\alpha})$ by $(z_1,\cdots,z_n).$
Let
$$\begin{array}{rl}
W=(w_1,\cdots,w_n)=&\sqrt{8}s(v_1(z),\cdots,v_n(z))
\sqrt{\frac{\cosh\theta -1}{\theta\sinh\theta}}\\
=&\sqrt{8}s\left((z_1,\cdots,z_n)A(z)^{*}+\ldots \right)
\sqrt{\frac{\cosh\theta -1}{\theta\sinh\theta}}.
\end{array}$$
It is easy to see that
$$\frac{\partial W}{\partial Z}=\frac{\partial (w_1,\cdots,w_n)}{\partial 
(z_1,\cdots,z_n)}$$ 
is non-degenerate. Thus 
$$\begin{array}{rl}
(s-\lim)\displaystyle\int_{N_{\epsilon}(p_{\alpha})}\hbox{str}&
\phi_0(\tau,t,q)dq=      
\lim_{t\rightarrow \infty}
\displaystyle\int_{N_{\epsilon}
(p_{\alpha})} \sqrt{\frac{1}{4\pi \tau^n}}\cdot \sqrt{\det \left(
\frac{\theta}{\sinh\theta}\right)}\\
&\quad \exp\{-\frac{W^2}{4\tau}\}
\exp\{s\sum v_{ij}(z)E_i^+ E_j^- \}(\det \frac{\partial W}{\partial Z})^{-1}
dW\\
=&\left(\hbox{str}\exp\{s\sum v_{ij}(z)E_i^+ E_j^- \}
(\det\frac{\partial W}{\partial Z})^{-1}
\cdot \sqrt{\det \left(\frac{\theta}{\sinh\theta}\right)} 
\right)_{z=0}.\end{array}$$
From
$$\frac{\partial W}{\partial Z}|_{z=0}=\left(2sA(z)^{*}
\sqrt{\frac{1}{\theta\sinh\theta}}\sqrt{2(\cosh\theta-1)}\right)_{z=0},$$
it follows that
$$\left(\det\frac{\partial W}{\partial Z}\right)^{-1}_{z=0}
=\left(
\sqrt{\det\frac{\sinh\theta}{\theta}}
\det(\sqrt{2(\cosh\theta-1)})^{-1}\right)_{z=0}.$$
So 
$$\begin{array}{rl}
(s-\lim)\displaystyle\int_{N_{\epsilon}(p_{\alpha})}\hbox{str}
\phi_0(\tau,t,q)dq=&
\left(\det \sqrt{2(\cosh\theta-1)})^{-1}\right)_{z=p_{\alpha}}\\
&\hbox{str}\exp\{s\sum v_{ij}(p_{\alpha})E_i^+ E_j^- \}.
\end{array}$$
Let us consider everything  when $s$ goes to zero. First we have
$$\begin{array}{rl}
[\det(\sqrt{2(\cosh\theta-1)})]_{z=0}
&=\det(\sqrt{\theta^2 +\ldots})_{z=0}\\
&=(2s)^{n}|\det(v_{ij}(p_{\alpha}))| +O(s^{n+1}).
\end{array}$$
Then let us consider the term
$\hbox{str}\exp\{s\sum v_{ij}(p_{\alpha})E_i^+ E_j^- \}.$
Note that 
$$ \hbox{str} ( E ^+ _{i_1} \cdots E ^+ _{i _\lambda } E ^- _{j _1} 
\cdots E ^- _{j _\mu } ) = 0, \quad \hbox{for} \ \lambda +\mu < 2n.$$
Hence 
$$\begin{array}{rl}
& \hbox{str} ( \exp (s \sum v _{jk} (p_{\alpha}) E ^+ _ j E ^- _ k ) ) 
= \frac{1}{n!} 
\hbox{str}\left( s\displaystyle\sum _{j,k} v _{jk}(p_{\alpha}) 
E ^+ _ j E ^- _k \right)^ n + O 
( s^{n +1}) \\
& = \frac{s^n }{n!} \sum \epsilon (j _1, \cdots, j _n) 
\epsilon (k _1, \cdots, k _n ) 
v _{j _1 k _1}(p_{\alpha})  
\cdots v _{j   _n k _n }(p_{\alpha})
\hbox{str} (E^+ _1 E ^- _1 
\cdots E ^+ _n E ^- _n ) \\
& \quad\quad + O (s^{n +1}), 
\end{array}$$
where the sum runs over all permutations $(j _1, \cdots, j 
_n),$ $(k_1, \cdots, k _n)$ of $(1,2, \cdots, n)$, and 
$\epsilon (j _1, \cdots, j 
_n)$ is equal to $1$ or $-1$ if the permutation $(j _1, \cdots, j _n )$ 
is even or odd, respectively. Therefore, by the Proposition 3' in [6]
$$
 \hbox{str} ( \exp (s \sum v _{jk} (p_{\alpha}) E ^+ _ j E ^- _ k ) 
 = ( 2s)^ n \cdot \det (v _{jk} (p)) + O (s^{n +1}). 
$$
From the above discussions we get
$$\begin{array}{rl}
\chi(M)=&
\lim_{s\rightarrow 0}(s-\lim) \displaystyle\int_M \hbox{str}
G(\tau,q,q,t)dq \\
=&\lim_{s\rightarrow 0}(s-\lim) \displaystyle\int_M \hbox{str}
\phi_0(\tau,t,q)dq\\
=&\displaystyle\sum_{p_{\alpha}\in \hbox{Zero}(V)}
\lim_{s \rightarrow 0}(s-\lim) \displaystyle\int_{N_{\epsilon} (p_{\alpha})} 
\hbox{str}\phi_0(\tau,t,q)dq \\
=&\displaystyle\sum_{p_{\alpha}\in \hbox{Zero}(V)}\lim_{s\rightarrow 0}
[\det\sqrt{2(\cosh\theta-1)}]^{-1}
\hbox{str}\exp\{s\sum v_{ij}(p_{\alpha})E_i^+ E_j^- \}\\
=&\displaystyle\sum_{p_{\alpha}\in \hbox{Zero}(V)}
\frac{\det (v_{ij}(p_{\alpha}))}{|\det (v_{ij}(p_{\alpha}))|}.
\end{array}$$
The theorem is proved.

\section{Appendix}

Let $M$ be an oriented Riemannian manifold of dimension $2$, $P:S(M)\rightarrow
M$ be the tangent sphere bundle of $M$. Three tangent vector fields
$X_1,X_2,X_3$ on $S(M)$ were well known for geometers (see the definitions
 in [4]). Let $\xi^t ,\eta^t$ be the
integral flows on $S(M)$, which correspond to $X_1,X_3$ respectively.
 A geodesic triangle can be described by a set of 
parameters $\{u,t,\theta,l,\gamma,b,\alpha \},$ where $u\in S(M),$ and
$\{t,\theta,l,\gamma,b,\alpha \}$ are arc lengths or angles,
such that the set of parameters satisfies
$$\xi^t \eta^{\pi -\alpha}\xi^b \eta^{\pi -\gamma}\xi^l \eta^{\pi -\theta}
u=u.
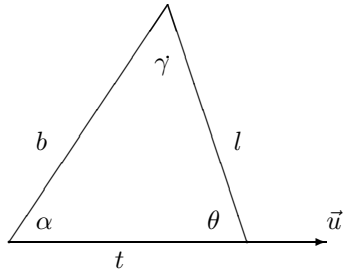
$$

\begin{picture}(300,100)(10,20)
\put(110,20){\line(1,0){90}}
\put (110,20){\line(2,3){60}}
\put (200,20){\line(-1,3){30}}
\put(200,20){\vector(1,0){30}}
\put(110,20){\circle*{.5}}
\put(200,20){\circle*{.5}}
\put(170,90){\circle*{.5}}
\put(150,10){\shortstack{$t$}}
\put(195,55){\shortstack{$l$}}
\put(120,55){\shortstack{$b$}}
\put(185,25){\shortstack{$\theta$}}
\put(165,85){\shortstack{$\gamma$}}
\put (120,25){\shortstack{$\alpha$}}
\put(230,25){$\vec{u}$}
\end{picture}
\vspace{.25in}

For a geodesic triangle, we denote
$$\begin{array}{rl}
B=P(u),&u_B=u,\\
C=P(\xi^l \eta^{\pi-\theta} u_{B}),& u_C =\xi^l \eta^{\pi-\theta }u_{B},\\
A=P(\xi^b \eta^{\pi-\gamma} u_{C}),& u_A =\xi^b \eta^{\pi-\gamma }u_{C},
\end{array}$$
and
$$\
[\vec{AB}]=\left(\begin{array}{cc}
                    1&0\\
                    0&H(|AB|,u_{B})\end{array}\right),$$
$$[\theta ] = \left(\begin{array}{ccc}
             -\cos \theta&-\sin \theta&0 \\
             \sin \theta&-\cos \theta&0 \\
              0&0&1  \end{array}\right),
$$
and so on, where $|AB|$ means the length of the geodesic arc $\overline{AB}$.
$H(t,u)$ is the unique solution to the following ODE
$$\left\{\begin{array}{rl}
 &\frac{d}{dt}X(t)=\left(\begin{array}{cc}
                            0&-k(p(t))\\
                            1&0
                            \end{array}\right)X(t)\\
 &X(0)=\left(\begin{array}{cc}
                            1&0\\
                            0&1
                            \end{array}\right)\end{array}\right. ,$$
where $X(t)$ are $2\times 2$ matrices, $u \in S(M), p(t)=P(\xi^{-t}u).$

{\bf Theorem(SAS trigonometry formulas)}\quad For not too big $t$ and $l$,
we can solve the geodesic triangle, i.e. there exist three functions
$$\begin{array}{rl}
\alpha =&\alpha(t,\theta,l,u),\\
\gamma =&\gamma(t,\theta,l,u),\\
b=&b(t,\theta,l,u),
\end{array}$$
such that
$$F(t,\theta,l,u)\equiv u,\qquad \forall (t,\theta,l,u),$$
where 
$$F(t,\theta,l,u)=
\xi^t \eta^{\pi -\alpha(t,\theta,l,u)}\xi^{b(t,\theta,l,u)}
 \eta^{\pi -\gamma(t,\theta,l,u)}\xi^l \eta^{\pi -\theta}u.$$
Then there hold
$$\begin{array}{rl}
(i)\quad &\left(\begin{array}{ccc}
          0&0&X_{1}(u)\alpha \\
          0&0&X_{2}(u)\alpha \\
          0&0&X_{3}(u)\alpha 
          \end{array}\right) [\vec{AB}]
-\left(\begin{array}{ccc}
          X_{1}(u)b&0&0 \\
          X_{2}(u)b&0&0 \\
          X_{3}(u)b&0&0
          \end{array}\right)[\alpha][\vec{AB}]\\
&\qquad \qquad +\left(\begin{array}{ccc}
          0&0&X_{1}(u)\gamma \\
          0&0&X_{2}(u)\gamma \\
          0&0&X_{3}(u)\gamma 
          \end{array}\right)[\vec{CA}][\alpha][\vec{AB}]\\
&=[\theta][\vec{BC}][\gamma][\vec{CA}][\alpha][\vec{AB}]
-\left(\begin{array}{ccc}

          1&0&0\\
          0&1&0 \\
          0&0&1
          \end{array}\right),\\
(ii)\quad &
\left(\begin{array}{ccc}
          0&0&\alpha_{t} \\
          0&0&\alpha_{\theta} \\
          0&0&\alpha_{l}
          \end{array}\right)[\vec{AB}]
-\left(\begin{array}{ccc}
          b_{t}&0&0\\
          b_{\theta}&0&0\\
          b_{l}&0&0
          \end{array}\right)[\alpha][\vec{AB}]\\
&\qquad \qquad +\left(\begin{array}{ccc}
          0&0&\gamma_{t} \\
          0&0&\gamma_{\theta} \\
          0&0&\gamma_{l}
          \end{array}\right)[\vec{CA}][\alpha][\vec{AB}]  \\
&=\left(\begin{array}{ccc}
          0&0&0 \\
          0&0&-1 \\
          0&0&0
          \end{array}\right)
[\vec{BC}][\gamma][\vec{CA}][\alpha][\vec{AB}]\\
&\qquad \qquad +\left(\begin{array}{ccc}
          0&0&0 \\
          0&0&0 \\
          1&0&0
          \end{array}\right)
[\gamma][\vec{CA}][\alpha][\vec{AB}]
 +\left(\begin{array}{ccc}
          1&0&0 \\
          0&0&0 \\
          0&0&0
          \end{array}\right),\end{array}$$
where $\alpha_{t}=\frac{\partial \alpha}{\partial t}.$

{\bf Remark}\quad If $M$ is a space form, then $\alpha,b,\gamma$
do not depend on $u$, hence (i) is just the last theorem in [4].

{\bf Proof}\quad Let 
$$\Phi:(0,t_{0}] \times S^{1} \times [0,l_{0}]\times S(M) \rightarrow
S(M)$$
be the projection. From $\Phi =F$ it follows 
$$\Phi_{*}=F_{*}.$$
From the trivial facts
$$\Phi_{*}(\frac{\partial}{\partial t})=
\Phi_{*}(\frac{\partial}{\partial \theta})=
\Phi_{*}(\frac{\partial}{\partial l})=0,\quad
\Phi_{*}(X_{i})=X_{i},$$
we get
$$F_{*}(\frac{\partial}{\partial t})=
F_{*}(\frac{\partial}{\partial \theta})=
F_{*}(\frac{\partial}{\partial l})=0,\quad
F_{*}(X_{i})=X_{i},$$
which turn out to be the desired (i) and (ii) due to concrete 
formulas of $F_{*}$.
For example,let us compute $F_{*}(\frac{\partial}{\partial l})$
$$\begin{array}{rl}
F_{*}(\frac{\partial}{\partial l})=&
(\xi^t \eta^{\pi -\alpha}\xi^b \eta^{\pi -\gamma}\xi^l \eta^{\pi -\theta}u)
_{*}(\frac{\partial}{\partial l})\\
=&(\xi^t \eta^{\pi -\alpha}\xi^b \eta^{\pi -\gamma})_{*}(X_{1}(\xi^l 
\eta^{\pi -\theta}u_{B}))\\
&+(\xi^t \eta^{\pi -\alpha}\xi^b)_{*}(X_{3}
( \eta^{\pi -\gamma}\xi^l \eta^{\pi -\theta}u_{B}))\cdot (-\gamma_{l})\\
&+(\xi^t \eta^{\pi -\alpha})_{*}(X_{1}
(\xi^b \eta^{\pi -\gamma}\xi^l \eta^{\pi -\theta}u_{B}))\cdot b_{l}\\
&+\xi^t_* (X_{3}( \eta^{\pi -\alpha}\xi^b \eta^{\pi -\gamma}
\xi^l \eta^{\pi -\theta}u_{B}))\cdot (-\alpha_{l}).
\end{array}$$
By using Lemma 7 in [4] we have 
$$\begin{array}{rl}
F_{*}(\frac{\partial}{\partial l})=
&(1,0,0)[\gamma][\vec{CA}][\alpha][\vec{AB}]\\
&-(0,0,\gamma_l)[\vec{CA}][\alpha][\vec{AB}]\\
&+(b_l ,0,0)[\alpha][\vec{AB}]\\
&-(0,0,\alpha_l)[\vec{AB}].
\end{array}$$
So $F_*(\frac{\partial}{\partial l})=0$ is just the last line of the 
equality (ii).

{\bf Theorem}\quad there hold
$$\begin{array}{rl}
(i)&b_l=\cos \gamma.\\
(ii)&b_{ll}=-\frac{H_{11}(b,u_A)}{H_{21}(b,u_A)}{\sin^{2}\gamma}.\\
(iii)&\hbox{there exist\quad} \epsilon>0 \hbox{\quad such that for } 0<t,l\leq \epsilon,\\
&(b^2)_{ll}\geq 1.
\end{array}$$

{\bf Proof}\quad The SAS formula shows
$$(0,0,\alpha_l)-(b_l,0,0)[\alpha]+(0,0,\gamma_l)[\vec{CA}][\alpha]
=(1,0,0)[\gamma][\vec{CA}][\alpha],$$
i.e.
$$\begin{array}{rl}
(0,0,&\alpha_l)-(b_l,0,0)\left(\begin{array}{ccc}
                                 -\cos \alpha&-\sin \alpha&0\\
                                 \sin \alpha&-\cos \alpha&0   \\
                                  0&0&1
                                \end{array}\right)\\
 &+(0,0,\gamma_l)\left(\begin{array}{cc}
                                 1&0   \\
                                  0&H(b,u_A)
                                \end{array}\right)
                \left(\begin{array}{ccc}
                                 -\cos \alpha&-\sin \alpha&0\\
                                 \sin \alpha&-\cos \alpha&0   \\
                                  0&0&1
                                \end{array}\right)        \\
=&(1,0,0)\left(\begin{array}{ccc}
                                 -\cos \gamma&-\sin \gamma&0\\
                                 \sin \gamma&-\cos \gamma&0   \\
                                  0&0&1
                                \end{array}\right)
         \left(\begin{array}{cc}
                    1&0\\
                    0&H(b,u_A)
               \end{array}\right)
         \left(\begin{array}{ccc}
                                 -\cos \alpha&-\sin \alpha&0\\
                                 \sin \alpha&-\cos \alpha&0   \\
                                  0&0&1
                                \end{array}\right),\end{array}$$
which is
$$\left\{ \begin{array}{l}
 b_l \cos\alpha  +H_{21}(b,u_A)\gamma_l \sin\alpha =\cos\gamma\cos\alpha
-H_{11}(b,u_A)
\sin\gamma\sin\alpha \\
b_l \sin\alpha  -H_{21}(b,u_A)\gamma_l\cos\alpha =\cos\gamma\sin\alpha
+H_{11}(b,u_A)\sin\gamma
\cos\alpha\\
\alpha_l +H_{22}(b,u_A)\gamma_l =-H_{12}(b,u_A)\sin\gamma.
\end{array}  \right.
$$ So  by the lemma 6 in [4]we get 
$$\left\{\begin{array}{rl}
b_l& =\cos\gamma\\
\gamma_l &=-\frac{H_{11}(b,u_A)}{H_{21}(b,u_A)}\sin\gamma\\
\alpha_l& =\frac{H_{22}H_{11}}
{H_{21}}\sin\gamma -H_{12}\sin\gamma\\
&=\frac{1}{H_{21}(b,u_A)}\sin\gamma
\end{array}\right.$$
And
$$\begin{array}{rl}
\frac{1}{2}(b^2)_{ll}=&(b_l)^2 +bb_{ll}={\cos^{2}\gamma}+(-sin\gamma \gamma_l)\\
=&{\cos^{2}\gamma} +\frac{b}{H_{21}}H_{11}{\sin^{2}\gamma}
=1+(\frac{b}{H_{21}}H_{11}-1){\sin^{2}\gamma}.
\end{array}$$
By the equation of $H(b,u_A)$ we know
$$\lim_{b\rightarrow 0}\frac{b}{H_{21}}=\lim_{b\rightarrow 0}H_{11}=1,$$
therefore there exists $\epsilon >0$, if $t,l \leq \epsilon,$
we have
$$b \leq t+l \leq 2 \epsilon,$$
so we can have
$$|\frac{b}{H_{21}} H_{11}-1|<\frac{1}{2}.$$
Thus
$$|(\frac{b}{H_{21}}H_{11} -1){\sin^{2}\gamma}|<\frac{1}{2},$$
$$\frac{1}{2}(b^2)_{ll}\geq 1- |(\frac{b}{H_{21}}H_{11} -1){\sin^{2}\gamma}|
\geq \frac{1}{2}.$$
The theorem is true now.

{\bf Theorem}\quad There exists an $\epsilon>0$, such that for any geodesic
$AB$ with $\rho(A,B)=|AB|=<\epsilon,$ for any $z\in M$, satisfying
$$\rho(z,A),\hbox{\quad} \rho(z,B)<\epsilon,$$
and for any $\lambda \in [0,1],$ we have
$$\mu \rho(z,A)^2+\lambda \rho(z,B)^2 -\lambda \mu \rho(A,B)^2 \geq
\frac{1}{4}\rho(o,z)^2,$$
where $\mu=1-\lambda,$ and $o$  is a point on the geodesic $AB$ such that
$$\frac{Ao}{AB}=\lambda.$$

{\bf Proof}\quad Let $\angle Aoz=\theta, l=\rho(o,z), \rho(A,B)=s,$
and let
$$f(l)=\lambda \rho(z,B)^2+\mu \rho(z,A)^2 -\lambda\mu l.$$
When $AB$ and $\theta$ are fixed,
$$\frac{d\rho(A,z)}{dl}|_{l=0}=\cos(\pi-\theta),\hbox{\quad} 
\frac{d\rho(B,z)}{dl}|_{l=0}=\cos\theta,$$
thus
$$f(0)=\lambda\mu^2 s^2 +\mu\lambda^2 s^2-\lambda\mu s^2 =0,$$and
$$\begin{array}{rl}
f'(0)&=\lambda\frac{d{\rho(z,B)}^2}{dl}|_{l=0}+\mu
\frac{d{\rho(z,A)}^2}{dl}|_{l=0}\\
&=\lambda\mu s \cos \theta+\mu\lambda\cos(\pi-\theta)=0.
\end{array}$$
Therefore
$$\begin{array}{rl}
f(l)&=\frac{1}{2}f''({\tilde l})l^2 =\frac{1}{2}\left(
\lambda \frac{d{\rho(z,B)}^2}{dl^2}
+\mu \frac{d{\rho(z,A)}^2}{dl^2} \right)|_{l={\tilde l}}l^2\\
&\geq \frac{1}{2}\left(\frac{\lambda}{2}+ \frac{\mu}{2}\right){\rho}^2
=\frac{\rho^2}{4}.
\end{array}$$
So the theorem is proved.

{\bf Acknowledgement.} The author is indebted to Prof 
H.Wu for supports to his works, and to Profs M.do Carmo,
A.Rigas, Macio A.F.R. Also, the author was partially 
supported by the the Science Fund of the
Chinese Academy of Sciences.

\begin{tabular}{l}
Institute of Mathematics\\[1mm]
Academia Sinica\\[1mm]
Beijing 100080\\[1mm]
PR China\\
email address: yuyl@math03.math.ac.cn
\end{tabular}

\end{document}